\documentclass[sn-mathphys-ay]{sn-jnl}

\usepackage{graphicx}%
\usepackage{multirow}%
\usepackage{amsmath,amssymb,amsfonts}%
\usepackage{amsthm}%
\usepackage{mathrsfs}%
\usepackage[title]{appendix}%
\usepackage{xcolor}%
\usepackage{textcomp}%
\usepackage{manyfoot}%
\usepackage{booktabs}%
\usepackage{algorithm}%
\usepackage{algorithmicx}%
\usepackage{algpseudocode}%
\usepackage{listings}%

\usepackage[capitalise,noabbrev]{cleveref} 
\usepackage{etoolbox} 
\usepackage{xstring} 
\usepackage{mathtools} 
\usepackage{subcaption} 
\usepackage{lmodern} 

\usepackage{tikz}
\usepackage{pgfplots}
\usepackage{pgfplotstable}
\usepgfplotslibrary{fillbetween}
\pgfplotsset{compat=1.18}
\immediate\write18{mkdir -p "externalize"}
\usepgfplotslibrary{external}
\tikzsetexternalprefix{externalize/}

\newcommand{\normal}[2]{\mathcal{N}\left(#1, #2\right)}
\newcommand{\GP}[2]{\mathcal{GP}\left(#1, #2\right)}
\newcommand{\linear}{\mathcal{L}}
\newcommand{\veclinear}{\boldsymbol{\linear}}
\newcommand{\information}{\mathcal{I}}
\newcommand{\vecinformation}{\boldsymbol{\information}}

\DeclareMathOperator*{\argmin}{arg\,min}

\newcommand{\boldtilde}[1]{\boldsymbol{\tilde{\text{$#1$}}}}
\newcommand{\boldhat}[1]{\boldsymbol{\hat{\text{$#1$}}}}
\newcommand{\boldstar}[1]{#1^*}

\newcommand{\matLambda}{\boldsymbol \Lambda}
\newcommand{\matSigma}{\boldsymbol \Sigma}
\newcommand{\matSigmaf}{\boldsymbol{\Sigma}_\mathbf{f}}
\newcommand{\matPhi}{\boldsymbol \Phi}

\newcommand{\vecsigma}{\boldsymbol \sigma}
\newcommand{\vecphi}{\boldsymbol \phi}
\newcommand{\vecpsi}{\boldsymbol \psi}

\newcommand{\matE}{\mathbf E}
\newcommand{\matF}{\mathbf F}
\newcommand{\matG}{\mathbf G}
\newcommand{\matGf}{\mathbf{G_f}}
\newcommand{\matH}{\mathbf H}
\newcommand{\matI}{\mathbf I}
\newcommand{\matK}{\mathbf K}
\newcommand{\matKc}{\mathbf{K_c}}
\newcommand{\matM}{\mathbf M}
\newcommand{\matP}{\mathbf P}
\newcommand{\matQ}{\mathbf Q}
\newcommand{\matX}{\mathbf X}
\newcommand{\vece}{\mathbf e}
\newcommand{\vecf}{\mathbf f}
\newcommand{\vecg}{\mathbf g}
\newcommand{\vecm}{\mathbf m}
\newcommand{\vecmf}{\mathbf{m_f}}
\newcommand{\vecu}{\mathbf u}
\newcommand{\vecuc}{\mathbf{u_c}}
\newcommand{\vecx}{\mathbf x}
\newcommand{\vecz}{\mathbf z}
\newcommand{\vecnull}{\mathbf 0}
\newcommand{\vecone}{\mathbf 1}

\newcommand{\subd}[1]{\ifblank{#1}{\mathbf{d}}{#1_\mathbf{d}}}
\newcommand{\subi}[1]{\ifblank{#1}{\mathbf{i}}{#1_\mathbf{i}}}
\newcommand{\subn}[1]{#1_\mathbf{n}}
\newcommand{\subii}[1]{\ifblank{#1}{\mathbf{ii}}{#1_\mathbf{ii}}}
\newcommand{\subid}[1]{#1_\mathbf{id}}
\newcommand{\subdi}[1]{#1_\mathbf{di}}
\newcommand{\subdd}[1]{#1_\mathbf{dd}}

\newcommand{\spaceV}{\mathcal{V}}
\newcommand{\spaceVh}{\mathcal{V}^h}
\newcommand{\spaceWh}{\mathcal{W}^h}
\newcommand{\spaceH}{H^1_0}
\newcommand{\funcuh}{u^h}
\newcommand{\funcvh}{v^h}

\newcommand{\diff}{\mathrm{d}}
\newcommand{\dx}{\, \diff x}
\newcommand{\vecdx}{\, \diff \vecx}

\newcommand{\vecdz}{\, \diff \vecz}

\definecolor{tudCyan}{HTML}{00A6D6}
\definecolor{tudWhite}{HTML}{FFFFFF}
\definecolor{tudBlack}{HTML}{000000}
\definecolor{tudNavy}{HTML}{0C2340}
\definecolor{tudSky}{HTML}{00B8C8}
\definecolor{tudBlue}{HTML}{0076C2}
\definecolor{tudPurple}{HTML}{6F1D77}
\definecolor{tudPink}{HTML}{EF60A3}
\definecolor{tudRed}{HTML}{A50034}
\definecolor{tudFire}{HTML}{E03C31}
\definecolor{tudOrange}{HTML}{EC6842}
\definecolor{tudYellow}{HTML}{FFB81C}
\definecolor{tudLime}{HTML}{6CC24A}
\definecolor{tudGreen}{HTML}{009B77}

\newcommand{\barplot}[5]{

    \begin{tikzpicture}
        \pgfplotstableset{
            create on use/x/.style={create col/copy column from table={img/1d-bar/M/mesh_Nf-64.dat}{x}}
        }

        \begin{axis}[
            xlabel={x},
            ylabel={u},
            ymin=-13,
            ymax=17,
            xtick={0.0,0.2,0.4,0.6,0.8,1.0},
            legend style={at={(1.0,1.05)},anchor=south east},
            legend cell align=left,
            legend columns=2,
            width=#3,
            height=#4,
        ]

            \addplot [name path=A, mark=none, opacity=0.0, forget plot] table [x=x, y expr=\thisrow{u_prior}-1.96*\thisrow{std_u_prior}] {img/1d-bar/#5/results_Nc-#1_Nf-64.dat};
            \addplot [name path=B, mark=none, opacity=0.0, forget plot] table [x=x, y expr=\thisrow{u_prior}+1.96*\thisrow{std_u_prior}] {img/1d-bar/#5/results_Nc-#1_Nf-64.dat};
            \addplot [name path=C, mark=none, opacity=0.0, forget plot] table [x=x, y expr=\thisrow{u_posterior}-1.96*\thisrow{std_u_posterior}] {img/1d-bar/#5/results_Nc-#1_Nf-64.dat};
            \addplot [name path=D, mark=none, opacity=0.0, forget plot] table [x=x, y expr=\thisrow{u_posterior}+1.96*\thisrow{std_u_posterior}] {img/1d-bar/#5/results_Nc-#1_Nf-64.dat};

            \addplot [color=tudOrange, opacity=0.25, forget plot] fill between [of=A and B];
            \addplot [color=tudBlue, opacity=0.25, forget plot] fill between [of=C and D];

            \foreach \i in {0,1,...,29} {\addplot [style=thin, opacity=0.25, color=tudOrange, forget plot] table [x=x, y=prior_sample_\i] {img/1d-bar/#5/samples-prior_Nc-#1_Nf-64.dat};}
            \foreach \i in {0,1,...,29} {\addplot [style=thin, opacity=0.25, color=tudBlue, forget plot] table [x=x, y=posterior_sample_\i] {img/1d-bar/#5/samples-posterior_Nc-#1_Nf-64.dat};}

            \addplot [name=Prior, style=very thick, mark=none, color=tudOrange] table [x=x, y=u_prior] {img/1d-bar/#5/results_Nc-#1_Nf-64.dat};
            \addplot [name=Posterior, style=very thick, mark=none, color=tudBlue] table [x=x, y=u_posterior] {img/1d-bar/#5/results_Nc-#1_Nf-64.dat};
            \addplot [name=Coarse solution, style=very thick, mark=none, color=tudRed] table [x=x, y=u_coarse] {img/1d-bar/#5/results_Nc-#1_Nf-64.dat};
            \addplot [name=Fine solution, style=very thick, mark=none, color=tudLime] table [x=x, y=u_fine] {img/1d-bar/#5/results_Nc-#1_Nf-64.dat};

            \addplot [name=Coarse nodes, only marks, color=black, domain=0:1, samples=#1+1, mark=*, mark size=1pt] {0};

            \IfEq{#2}{1}{
                \legend{Prior, Posterior, Coarse solution, Fine solution, Coarse nodes}
            }{}
        \end{axis}
    \end{tikzpicture}
}

\begin{document}

\title{A Bayesian Approach to Modeling Finite Element Discretization Error}


\author*[1]{\fnm{Anne} \sur{Poot}}\email{a.poot-1@tudelft.nl}
\author[2]{\fnm{Pierre} \sur{Kerfriden}}\email{pierre.kerfriden@minesparis.psl.eu}
\author[1]{\fnm{Iuri} \sur{Rocha}}\email{i.rocha@tudelft.nl}
\author[1]{\fnm{Frans} \spfx{van der} \sur{Meer}}\email{f.p.vandermeer@tudelft.nl}

\affil*[1]{\orgdiv{Faculty of Civil Engineering and Geosciences}, \orgname{Delft University of Technology}, \orgaddress{\street{Stevinweg 1}, \city{Delft}, \postcode{2628 CN}, \country{the Netherlands}}}
\affil[2]{\orgdiv{Centres des Matériaux}, \orgname{Mines Paris -- PSL}, \orgaddress{\street{63-65 Rue Henri Auguste Desbruères}, \city{Évry}, \postcode{91100}, \country{France}}}

\abstract{
    In this work, the uncertainty associated with the finite element discretization error is modeled following the Bayesian paradigm.
    First, a continuous formulation is derived, where a Gaussian process prior over the solution space is updated based on observations from a finite element discretization.
    To avoid the computation of intractable integrals, a second, finer, discretization is introduced that is assumed sufficiently dense to represent the true solution field.
    A prior distribution is assumed over the fine discretization, which is then updated based on observations from the coarse discretization.
    This yields a posterior distribution with a mean that serves as an estimate of the solution, and a covariance that models the uncertainty associated with this estimate.
    Two particular choices of prior are investigated: a prior defined implicitly by assigning a white noise distribution to the right-hand side term, and a prior whose covariance function is equal to the Green's function of the partial differential equation.
    The former yields a posterior distribution with a mean close to the reference solution, but a covariance that contains little information regarding the finite element discretization error.
    The latter, on the other hand, yields posterior distribution with a mean equal to the coarse finite element solution, and a covariance with a close connection to the discretization error.
    For both choices of prior a contradiction arises, since the discretization error depends on the right-hand side term, but the posterior covariance does not.
    We demonstrate how, by rescaling the eigenvalues of the posterior covariance, this independence can be avoided.
}


\maketitle

\section{Introduction}\label{sec:introduction}
In recent years, the Bayesian paradigm has become a popular framework to perform uncertainty quantification.
It has found its application in global optimization~\citep{mockus_bayesian_1989}, inverse modeling~\citep{stuart_inverse_2010} and data assimilation~\citep{law_data_2015} contexts, among others.
Commonly, given some numerical model, a prior distribution is assumed over its parameters, and the Bayesian paradigm provides a consistent framework to estimate these parameters and to quantify and propagate their associated uncertainty.
It should be noted, however, that even if complete certainty could be obtained over the model parameters, there would still be a remaining uncertainty to the solution due to approximations made in the numerical model.
This key observation is what underpins the current trend towards probabilistic numerics.

At the core of probabilistic numerics, the estimation of an unknown field is recast as a statistical inference problem, which allows for the estimation of the field with some uncertainty measure~\citep{larkin_gaussian_1972,Diaconis1988}.
Early examples of the application of Bayesian probabilistic numerics include computing integrals~\citep{ohagan_bayeshermite_1991} and solving ordinary differential equations~\citep{skilling_bayesian_1992}.
More recently, following a ``call to arms'' from \citet{hennig_probabilistic_2015-1}, a large push has been made to apply this framework to a wide range of problems, ranging from solving linear systems~\citep{hennig_probabilistic_2015,cockayne_bayesian_2019-1,wenger_probabilistic_2020} to quadrature~\citep{karvonen_classical_2017,briol_probabilistic_2017} to solving ordinary differential equations~\citep{schober_probabilistic_2014,hennig_probabilistic_2014,teymur_probabilistic_2016}.
For a general overview of the current state-of-the-art of probabilistic numerics, the reader is referred to \citet{hennig_probabilistic_2022}.
Most relevant for the work presented in this paper are the probabilistic numerical methods that have been developed for the solving of partial differential equations, which can be roughly divided into two categories: meshfree probabilistic solvers, and solver-perturbing error estimators.

The first category~\citep{chkrebtii_bayesian_2016,cockayne_probabilistic_2017,raissi_numerical_2018,wang_bayesian_2021} can be seen as a way to find solutions to partial differential equations directly from the strong form in a Bayesian manner.
A prior is assumed over the solution field, which is updated by evaluating its derivatives on a grid of collocation points, allowing for a solution to be obtained without needing to apply a finite element discretization over the domain.
This approach to solving partial differential equations shares some similarities with Bayesian physics-informed neural networks~\citep{raissi_physics-informed_2019,yang_b-pinns_2021}, the main difference lying in the function that is being fitted at the collocation points.
The way in which these meshfree solvers relate to traditional collocation methods is similar to the way in which Bayesian physics-informed neural networks relate to their deterministic counterparts.

The second category~\citep{conrad_statistical_2017,kersting_active_2018,lie_strong_2019} is focused on estimating the discretization error of traditional solvers for differential equations.
For ordinary differential equations, the usual time integration step is taken, after which the solution is perturbed by adding Gaussian noise, representing the uncertainty in the time integration result.
Similarly, for partial differential equations, the traditional spatial discretization is perturbed using small support Gaussian random fields, which reflect the uncertainty introduced by the mesh.
In~\citet{abdulle_random_2020,abdulle_probabilistic_2021}, a similar approach is taken, but rather than adding noise to the solution, an uncertainty is introduced by perturbing the time step size or finite element discretization.
A more formal mathematical basis for probabilistic numerical methods can be found in~\citet{cockayne_bayesian_2019}, where a more rigorous definition of the term is outlined and a common framework underpinning these two seemingly separate categories is established.

It is worth noting that these probabilistic numerical methods are a deviation from traditional error estimators~\citep{babuska_analysis_1979,babuska_feedback_1987,zienkiewicz_simple_1987}, as they embed the model error into the method itself, rather than estimate it a posteriori.
This inherently affects the model output, which depending on the context can be a desirable or undesirable property.
In~\citet{rouse_probabilistic_2021}, a method is presented to obtain full-field error estimates by assuming a Gaussian process prior over the discretization error, and updating it based on a set of traditional estimators of error in quantities of interest.
This way, a distribution representing the finite element discretization error can be obtained in a non-intrusive manner.

The shared goal of these methods is to accurately describe the errors made due to limitations of our numerical models, though their method of modeling error differs.
At the core, the meshfree probabilistic solvers model error as the result of using a finite number of observations to obtain a solution to an infinite-dimensional problem.
The solver-perturbing error estimators, on the other hand, take an existing discretization, like the one used in the finite element method, and assign some uncertainty measure to the existing solver.
This begs the question:
what happens if the methodology from the meshfree probabilistic solvers is applied to existing mesh-based solvers of partial differential equations?
Little research has thus far been conducted to answer this question, though two particular works are worth pointing out.

A brief remark is made in~\citet{bilionis_probabilistic_2016} describing a Bayesian probabilistic numerical method whose posterior mean is equivalent to the finite element solution.
However, this idea is then discarded due to infinite variances arising in the posterior distribution.
In~\citet{pfortner_physics-informed_2023}, the probabilistic meshfree solvers from \citet{cockayne_probabilistic_2017} are generalized to methods of weighted residuals, which includes the finite element method.
Of particular relevance to our work is their construction of prior distributions whose posterior mean is guaranteed to be equivalent to the usual finite element solution.
Doing this would allow one to replace the traditional finite element solver with the probabilistic one, in order to quantify the finite element discretization error.
However, their experimental results are limited to one-dimensional test cases, possibly because the application of their formulation to unstructured triangular or quadrilateral meshes would result in integrals in the information operator that are computationally too expensive.

In this work, we propose a probabilistic numerical method for the modeling of finite element discretization error.
The solution is endowed with a Gaussian process prior, which is then updated based on observations of the right-hand side from a finite element discretization.
This allows for the approximation of the true solution while including the uncertainty resulting from the finite discretization that is applied.
Rather than work directly with the Gaussian process distribution over the exact solution space, we introduce a second discretization over the domain that is fine enough to represent the exact solution.
This second discretization helps to avoid the infinite variances brought up in~\citet{bilionis_probabilistic_2016} as well as the computationally expensive integrals from~\citet{pfortner_physics-informed_2023}.
We present a class of priors that naturally accounts for the smoothness of the partial differential equation at hand, and show how the assembly of large full covariance matrices can be avoided.
A particular focus of this work is on the relationship between the posterior covariance of our formulation and the finite element discretization error.
The relationship between these two quantities is often left to intuition, reasoning along the lines that since the posterior covariance contains remaining model uncertainty, it must reflect the discretization error.
We challenge this assumption and investigate more thoroughly which conditions need to be met before the posterior covariance can reasonably be said to capture the finite element discretization error.

The underlying goal of the development of a Bayesian model for the finite element discretization error is to enable the propagation of discretization error to quantities of interest through the computational pipelines that arise in multiscale modeling, inverse modeling and data assimilation settings.
This consistent treatment of discretization error in turn allows for more informed decisions to be made about its impact on the model output.
To give a concrete example, in~\citet{girolami_statistical_2021}, a Bayesian framework for the assimilation of measurement data and finite element models is presented.
Within this framework, a model misspecification component is defined, which is endowed with a squared-exponential Gaussian process prior.
The Bayesian formulation of the finite element method that we derive in this work would allow for a more informative choice of prior distribution over the model misspecification component, for example by separating out the discretization error from the error associated with other modeling assumptions.

The outline of this paper is as follows:
in \cref{sec:bfem}, we derive our Bayesian formulation of the finite element method.
This is followed by a discussion on the choice of prior covariance in \cref{sec:choice-of-prior}, where two different choices of prior distribution are investigated.
Two examples, a one-dimensional tapered bar and a two-dimensional perforated plate, are showcased throughout this section to validate the conclusions drawn from theory.
Finally, in \cref{sec:conclusions}, the conclusions of this paper are drawn and discussed.

\section{Bayesian Finite Element Method}\label{sec:bfem}
In this section, the proposed Bayesian version of the finite element method is derived.
Although the method is applicable to a broad range of linear elliptic partial differential equations, for the purposes of demonstration, we will consider Poisson's equation:
\begin{equation}\label{eq:strong-form}
\begin{aligned}
    -\Delta u(\vecx) &= f(\vecx) & &\text{ in } \Omega\\
    u(\vecx) &= 0 & &\text{ on } \partial \Omega
\end{aligned}
\end{equation}
Here, $\Omega$ and $\partial \Omega$ are the domain and its boundary, respectively.
$u(\vecx)$ and $f(\vecx)$ are the solution and forcing term, which are linked through the Laplace operator $\Delta$.

\subsection{Continuous formulation}\label{subsec:continuous-form}
We will start with the derivation of a continuous posterior distribution over the solution space conditioned on the finite element force vector, largely following \citet{bilionis_probabilistic_2016}.
As usual, the problem is restated in its weak formulation:
\begin{equation}\label{eq:weak-form}
\begin{aligned}
    \int_{\Omega} \nabla u(\vecx) \cdot \nabla v(\vecx) \vecdx &= \int_{\Omega} f(\vecx) v(\vecx) \vecdx & &\forall v(\vecx) \in \spaceV
\end{aligned}
\end{equation}
We search $u(\vecx) \in \spaceV$, where $\spaceV=\spaceH$ is a Sobolev space of functions over $\Omega$ that are weakly once-differentiable and vanish at the boundary $\partial \Omega$.
This space is equipped with an inner product and thus also forms a Hilbert space.
Now, a discretization is defined over the domain using a set of locally supported shape functions $\{\psi_i(\vecx)\}_{i=1}^m$, which span a finite-dimensional space $\spaceWh \subset \spaceV$.
The test function $\funcvh(\vecx)$ can be defined in terms of these shape functions:
\begin{equation}\label{eq:test-functions}
\begin{aligned}
    \funcvh(\vecx) = \sum_{i=1}^m v_i \psi_i(\vecx) &&\text{ with } \psi_i(\vecx) \in \spaceWh
\end{aligned}
\end{equation}
Since \cref{eq:weak-form} has to hold for all $\funcvh(\vecx) \in \spaceWh$, the weights $v_i$ can be chosen at will.
Substituting \cref{eq:test-functions} into \cref{eq:weak-form}, a finite set of $m$ equations in constructed by choosing $v_i = \delta_{ij}$ for the $j$th equation, where $\delta_{ij}$ is the Kronecker delta function.
This yields the entries of the finite element force vector $\vecg$:
\begin{equation}\label{eq:coarse-force-vector}
\begin{aligned}
    g_i = \int_{\Omega} f(\vecx) \psi_i(\vecx) \vecdx
\end{aligned}
\end{equation}
We can relate the solution $u(\vecx)$ to the force vector $\vecg$ via the linear operator $\veclinear$:
\begin{equation}\label{eq:linear-relationship}
\begin{aligned}
    \veclinear\left[u(\vecx)\right] = \vecg
\end{aligned}
\end{equation}
where $\veclinear\left[u(\vecx)\right] = \begin{bsmallmatrix} \linear_1\left[u(\vecx)\right] & \linear_2\left[u(\vecx)\right] & \dots & \linear_m\left[u(\vecx)\right] \end{bsmallmatrix}^T$ is given by
\begin{equation}\label{eq:linear-observation-operator}
\begin{aligned}
    \linear_i [u(\vecx)] = \int_{\Omega} \nabla u(\vecx) \cdot \nabla \psi_i(\vecx) \vecdx
\end{aligned}
\end{equation}

A centered Gaussian process with a positive definite covariance function $k(\vecx, \vecx')$ is now assumed over the solution $u(\vecx)$:
\begin{equation}\label{eq:prior-distribution-continuous}
\begin{aligned}
    u(\vecx) \sim \GP{0}{k(\vecx, \vecx')}
\end{aligned}
\end{equation}
Because we have a linear map $\veclinear$ from $u(\vecx)$ to $\vecg$, conditioning $u(\vecx)$ on $\vecg$ yields another Gaussian process distribution~\citep{pfortner_physics-informed_2023}:
\begin{equation}\label{eq:posterior-distribution-continuous}
\begin{aligned}
    u(\vecx) | \vecg \sim \GP{m^*(\vecx)}{k^*(\vecx,\vecx')}
\end{aligned}
\end{equation}
Here, the posterior mean function $m^*(\vecx)$ and covariance function $k^*(\vecx,\vecx')$ are given by\footnote{To avoid confusion when applying $\veclinear$ to the covariance function $k(\vecz, \vecz')$, we use $\veclinear$ and $\veclinear'$ to denote that gradients and integrals are computed with respect to $\vecz$ and $\vecz'$, respectively.}:
\begin{equation}\label{eq:posterior-moments-continuous}
\begin{aligned}
    m^*(\vecx) &= \veclinear' \left[k(\vecx, \vecz')\right] \left(\veclinear\left[\veclinear'\left[k(\vecz, \vecz')\right]\right]\right)^{-1} \vecg \\
    k^*(\vecx,\vecx') &= k(\vecx, \vecx') - \veclinear'\left[k(\vecx, \vecz')\right] \left(\veclinear\left[\veclinear'\left[k(\vecz, \vecz')\right]\right]\right)^{-1} \veclinear\left[k(\vecz, \vecx')\right]
\end{aligned}
\end{equation}
The posterior mean function $m^*(\vecx)$ provides an full-field estimate of the solution $u(\vecx)$.
The posterior covariance function $k^*(\vecx, \vecx')$ indicates the uncertainty associated with this estimate due to the fact that it was obtained using only a finite set of shape functions.
Since the finite discretization is the only source of uncertainty in our model, we can intuit some association between this posterior covariance and the finite element discretization error.

The formulation presented thus far can be contextualized in the method of weighted residuals framework presented in \citet{pfortner_physics-informed_2023}.
Specifically, our continuous formulation is equivalent to choosing the information operator $\vecinformation \left[u(\vecx)\right] = \begin{bsmallmatrix} \information_1\left[u(\vecx)\right] & \information_2\left[u(\vecx)\right] & \dots & \information_m\left[u(\vecx)\right] \end{bsmallmatrix}^T$ in their framework to be given by:
\begin{equation}\label{eq:information-operator}
\begin{aligned}
    \information_i\left[u(\vecx)\right] = \int_\Omega \nabla u(\vecx) \cdot \nabla \psi_i(\vecx) \vecdx - \int_\Omega f(\vecx) \psi_i(\vecx) \vecdx
\end{aligned}
\end{equation}

Unfortunately, the integrals that arise in the expressions for the posterior mean and covariance functions in \cref{eq:posterior-moments-continuous} are generally intractable.
For some arbitrary covariance function $k(\vecx, \vecx')$, the integration over the shape functions $\psi_i(\vecx)$ and $\psi_j(\vecx')$ cannot be performed without putting severe restrictions on which shape functions are permitted.
This in turn puts severe constraints on the domain shape, which undercuts the core strength of the finite element method, namely its ability to solve partial differential equations on complicated domains.
On the other hand, we can design the covariance function such that these integrals do become tractable, for example by following~\citet{bilionis_probabilistic_2016} and setting $k(\vecx, \vecx') = G(\vecx, \vecx')$, or following~\citet{owhadi_bayesian_2015} and setting $k(\vecx, \vecx') = \int_\Omega \int_\Omega G(\vecx, \vecz) G(\vecx', \vecz') \delta(\vecz - \vecz') \vecdz \vecdz'$, where $\delta(\vecx)$ is a Dirac delta function.
However, in both of these expressions, the Green's function $G(\vecx, \vecx')$ associated with the operator $-\Delta$ is required, which is generally not available for a given partial differential equation.
Since our aim is to develop a general Bayesian framework for modeling finite element discretization error, a new approach is needed that does not impose restrictions on the choice of shape functions or require access to the Green's function.

\subsection{Discretized formulation}\label{subsec:petrov-galerkin}
This motivates us to approximate $u(\vecx)$ in the finite-dimensional space $\spaceVh$ spanned by a second set of locally supported shape functions $\{\phi_j(\vecx)\}_{j=1}^n$.
This defines the trial function $\funcuh(\vecx)$:
\begin{equation}\label{eq:trial-functions}
\begin{aligned}
    u(\vecx) \approx \funcuh(\vecx) = \sum_{j=1}^n u_j \phi_j(\vecx) &&\text{ with } \phi_j(\vecx) \in \spaceVh
\end{aligned}
\end{equation}
Note that this is not the same set of shape functions as the one used to define the force vector in \cref{eq:coarse-force-vector}.
In fact, since our aim is to model the discretization error that arises by choosing $v(\vecx) \in \spaceWh$ rather than and $v(\vecx) \in \spaceV$, it is important that the error associated with the projection of an arbitrary function $w(\vecx) \in \spaceV$ onto $\spaceVh$ is small compared to the error associated with its projection onto $\spaceWh$.
Loosely speaking, we assume that $\spaceVh$ is sufficiently expressive to serve as a stand-in for $\spaceV$.

Substituting \cref{eq:test-functions,eq:trial-functions} into \cref{eq:weak-form} yields the matrix formulation of the problem:
\begin{equation}\label{eq:matrix-form}
\begin{aligned}
    \matH \vecu &= \vecg
\end{aligned}
\end{equation}
The elements of the stiffness matrix $\matH$ are given by:
\begin{equation}\label{eq:stiffness-matrix}
\begin{aligned}
    H_{ij} &= \int_{\Omega} \nabla \psi_i(\vecx) \cdot \nabla \phi_j(\vecx) \vecdx \\
\end{aligned}
\end{equation}
The assumption that $\spaceVh$ is more expressive than $\spaceWh$ implies that $\vecu$ will have a larger dimensionality than $\vecg$ and thus that $\matH$ is a rectangular matrix and that \cref{eq:matrix-form} describes an underdetermined system.
However, the fact that this system of equations has an infinite set of solutions need not pose a problem, due to the regularizing effect of the prior assumed over $u(\vecx)$.

Since the solution field $u(\vecx)$ has been reduced from the infinite-dimensional space $\spaceV$ to the finite-dimensional $\spaceVh$, the distribution assumed over the solution in \cref{eq:prior-distribution-continuous} needs to be reduced accordingly.
Instead of an infinite-dimensional Gaussian process, we obtain a finite-dimensional zero-mean normal distribution with a positive definite covariance matrix $\matSigma$:
\begin{equation}\label{eq:prior-distribution-finite}
\begin{aligned}
    \vecu \sim \normal{0}{\matSigma}
\end{aligned}
\end{equation}
The joint distribution of $\vecu$ and $\vecg$ is now given by:
\begin{equation}\label{eq:base-joint-distribution}
\begin{aligned}
    \begin{bmatrix}
        \vecg \\
        \vecu
    \end{bmatrix}
    =
    \begin{bmatrix}
        \matH \vecu \\
        \vecu
    \end{bmatrix}
    \sim \normal{\vecnull}{
        \begin{bmatrix}
            \matH \matSigma \matH^T & \matH \matSigma \\
            \matSigma \matH^T & \matSigma
        \end{bmatrix}
    }
\end{aligned}
\end{equation}
Conditioning $\vecu$ on $\vecg$ yields the following posterior distribution:
\begin{equation}\label{eq:base-posterior-distribution}
\begin{aligned}
    \vecu | \vecg \sim \normal{\boldstar{\vecm}}{\boldstar{\matSigma}}
\end{aligned}
\end{equation}
Here, the posterior mean vector $\boldstar{\vecm}$ and covariance matrix $\boldstar{\matSigma}$ are given by:
\begin{equation}\label{eq:base-posterior-mean}
\begin{aligned}
    \boldstar{\vecm} &= \matSigma \matH^T \left(\matH \matSigma \matH^T\right)^{-1} \vecg \\
    \boldstar{\matSigma} &= \matSigma - \matSigma \matH^T \left(\matH \matSigma \matH^T \right)^{-1} \matH \matSigma
\end{aligned}
\end{equation}

Similar to the continuous formulation presented in section~\ref{subsec:continuous-form}, $\boldstar{\vecm}$ can be interpreted as providing an estimate of the solution $u(\vecx)$ in the fine space $\spaceVh$, while observing the right-hand side $f(\vecx)$ only in the coarse space $\spaceWh$.
The posterior covariance matrix $\boldstar{\matSigma}$ then provides an indication of the uncertainty associated with this estimate due to the fact that only observations from the coarse mesh are used to obtain this estimate.
Note that if the test and trial spaces are chosen to be the same (i.e.\ $\spaceWh = \spaceVh$), $\boldstar{\matSigma}$ reduces to a null matrix, reflecting the fact that there no longer exists a discretization error between $\spaceVh$ and $\spaceWh$.

\subsection{Hierarchical shape functions}\label{subsec:hierarchical-shape-functions}
Thus far, the only requirement that has been put on the choice of $\spaceVh$ and $\spaceWh$ is that the error between $\spaceV$ and $\spaceVh$ is small compared to the error between $\spaceV$ and $\spaceWh$.
We now add a second restriction, namely that $\spaceWh \subset \spaceVh$.
This defines a hierarchy between these two spaces, and implies that any function defined in $\spaceWh$ can be expressed in $\spaceVh$.
One way to ensure this hierarchy in practice is to first define a coarse mesh corresponding to $\spaceWh$, and then refine it hierarchically to obtain a fine mesh corresponding to $\spaceVh$.
Alternatively, it is possible to use only a single mesh, and use linear and quadratic shape functions over the same finite elements to define $\spaceWh$ and $\spaceVh$, respectively.

From the hierarchy between $\spaceVh$ and $\spaceWh$, it follows that the basis functions that span the coarse space $\spaceWh$ can be written as a linear combination of the basis functions that span the fine space $\spaceVh$.
In other words, there exists a matrix\footnote{Note that $\matPhi$ has been defined in terms of its transpose in order to make expressions in later sections consistent with common notation for least squares, proper orthogonal decomposition, and so on.} $\matPhi^T$ that maps a vector of fine shape functions $\vecphi(\vecx) = \begin{bsmallmatrix} \phi_1(\vecx) & \phi_2(\vecx) & \dots & \phi_n(\vecx) \end{bsmallmatrix}^T$ to a vector of coarse shape functions $\vecpsi(\vecx) = \begin{bsmallmatrix} \psi_1(\vecx) & \psi_2(\vecx) & \dots & \psi_m(\vecx) \end{bsmallmatrix}^T$:
\begin{equation}\label{eq:shape-function-hierarchy}
\begin{aligned}
    \vecpsi(\vecx) = \matPhi^T \vecphi(\vecx)
\end{aligned}
\end{equation}
This allows \cref{eq:stiffness-matrix} to be rewritten as:
\begin{equation}\label{eq:stiffness-matrix-hierarchical}
\begin{aligned}
    H_{ij} &= \int_{\Omega} \nabla \sum_{k=1}^n \Phi_{ki} \phi_k(\vecx) \cdot \nabla \phi_j(\vecx) \vecdx \\
    &= \sum_{k=1}^n \Phi_{ki} \int_{\Omega} \nabla \phi_k(\vecx) \cdot \nabla \phi_j(\vecx) \vecdx
\end{aligned}
\end{equation}
As a result, $\matH$ can be expressed as:
\begin{equation}\label{eq:phi-H-K-relation}
\begin{aligned}
    \matH = \matPhi^T \matK
\end{aligned}
\end{equation}
where $\matK$ is the fine-scale (square and symmetric) stiffness matrix that would follow if both trial and test functions came from the fine space $\spaceVh$:
\begin{equation}\label{eq:bg-stiffness-matrix}
\begin{aligned}
    K_{ij} &= \int_{\Omega} \nabla \phi_i(\vecx) \cdot \nabla \phi_j(\vecx) \vecdx \\
\end{aligned}
\end{equation}
Following a similar line of reasoning, the coarse stiffness matrix $\matKc$, that would be found if both trial and test functions came from the coarse space $\spaceWh$, can be written in terms of $\matPhi$ and $\matK$:
\begin{equation}\label{eq:phi-Kc-K-relation}
\begin{aligned}
    \matKc &= \matPhi^T \matK \matPhi
\end{aligned}
\end{equation}

Similarly to \cref{eq:stiffness-matrix-hierarchical}, we can rewrite \cref{eq:coarse-force-vector} as:
\begin{equation}\label{eq:force-vector-hierarchical}
\begin{aligned}
    g_{i} &= \int_{\Omega} f(\vecx) \sum_{k=1}^n \Phi_{ki} \phi_k(\vecx) \vecdx \\
    &= \sum_{k=1}^n \Phi_{ki} \int_{\Omega} f(\vecx) \phi_k(\vecx) \vecdx
\end{aligned}
\end{equation}
And so, $\vecg$ can be expressed as:
\begin{equation}\label{eq:phi-g-f-relation}
\begin{aligned}
    \vecg = \matPhi^T \vecf
\end{aligned}
\end{equation}
where $\vecf$ is the fine-scale force vector that arises by integrating the forcing term over the fine-scale test functions:
\begin{equation}\label{eq:fine-force-vector}
\begin{aligned}
    f_i &= \int_{\Omega} f(\vecx) \phi_i(\vecx) \vecdx
\end{aligned}
\end{equation}

Finally, we define the reference solution $\boldhat{\vecu}$ as the solution to the fine-scale system of equations that is obtained by choosing both the test and trial spaces to be the fine space $\spaceVh$:
\begin{equation}\label{eq:reference-solution}
\begin{aligned}
    \matK \boldhat{\vecu} = \vecf
\end{aligned}
\end{equation}
In the remainder of this work, discretization error is defined with respect to $\boldhat{\vecu}$.
Specifically, the finite element discretization error $\vece$ is defined as the difference between the fine-scale reference solution $\boldhat{\vecu}$ and the coarse-scale solution projected to the fine space:
\begin{equation}\label{eq:discretization-error-definition}
\begin{aligned}
    \vece = \matK^{-1} \vecf - \matPhi \matKc^{-1} \vecg
\end{aligned}
\end{equation}

\subsection{Boundary conditions}\label{subsec:dirichlet-bcs}
It is worth considering how the application of boundary conditions in the fine space translates to the shape functions in the coarse space.
To do this, $\vecphi(\vecx)$ is split into $\subi{\vecphi}(\vecx)$ and $\subd{\vecphi}(\vecx)$, where the subscript $\subd{}$ refers to the nodes on the part of the boundary where Dirichlet conditions are applied, and the subscript $\subi{}$ refers to all other nodes (i.e.\ both internal nodes and non-Dirichlet boundary nodes).
This could be considered abuse of notation, since $\spaceVh \subset \spaceV$, which is already constrained by the Dirichlet boundary conditions, so from this point of view, $\subd{\vecphi}(\vecx)$ should not exist.
However, in most practical finite element implementations, shape functions are assigned to the boundary nodes as well in order to facilitate the inclusion of inhomogeneous boundary conditions in the model.

The boundary conditions in the coarse space follow from $\subd{\vecphi}(\vecx)$ and $\matPhi$, since $\subd{\vecpsi}(\vecx)$ is defined as the elements of $\vecpsi(\vecx)$ where the rows of $\matPhi$ belonging to $\subd{\vecphi}(\vecx)$ have non-zero entries.
As a result, \cref{eq:shape-function-hierarchy} can be split as follows:
\begin{equation}\label{eq:phi-split}
\begin{aligned}
    \begin{bmatrix}
        \subi{\vecpsi}(\vecx) \\
        \subd{\vecpsi}(\vecx)
    \end{bmatrix}
    =
    \begin{bmatrix}
        \subii{\matPhi}^T & \vecnull \\
        \subid{\matPhi}^T & \subdd{\matPhi}^T
    \end{bmatrix}
    \begin{bmatrix}
        \subi{\vecphi}(\vecx) \\
        \subd{\vecphi}(\vecx)
    \end{bmatrix}
\end{aligned}
\end{equation}
Note that the fact that $\subdi{\matPhi} = \vecnull$ does not introduce any loss of generality:
any non-zero element of $\subdi{\matPhi}$ would by definition of $\subd{\vecpsi}(\vecx)$ be an element of $\subdd{\matPhi}$, not $\subdi{\matPhi}$.
From \cref{eq:phi-H-K-relation,eq:phi-split}, it follows that:
\begin{equation}\label{eq:int-phi-H-K-relation}
\begin{aligned}
    \subii{\matH} &= \subii{\matPhi}^T \subii{\matK}
\end{aligned}
\end{equation}
Similarly, from \cref{eq:phi-g-f-relation,eq:phi-split}, it follows that:
\begin{equation}\label{eq:int-phi-g-f-relation}
\begin{aligned}
    \subi{\vecg} &= \subii{\matPhi}^T \subi{\vecf}
\end{aligned}
\end{equation}

Commonly, Dirichlet boundary conditions are enforced by eliminating the corresponding degrees of freedom, and solving the system that remains.
Due to the simple relation that $\subii{\matPhi}$ provides between $\subii{\matH}$ and $\subii{\matK}$ (\cref{eq:int-phi-H-K-relation}) as well as $\subi{\vecg}$ and $\subi{\vecf}$ (\cref{eq:int-phi-g-f-relation}), all relationships described in \cref{subsec:petrov-galerkin,subsec:hierarchical-shape-functions} still hold when applied only to the internal nodes of the system.
From this point onward, we will therefore only consider the internal nodes of the system.
This also means that only the part of the covariance matrix related to the internal nodes $\subii{\matSigma}$ needs to be considered, and so the requirement of positive definiteness of $\matSigma$ can be relaxed to a requirement of positive definiteness of only $\subii{\matSigma}$.
The subscripts $\subi{}$ (for vectors) and $\subii{}$ (for matrices) will be left implied in order to declutter the notation.

In the remainder of this paper, we will limit ourselves to partial differential equations with homogeneous boundary conditions.
However, the method can easily be extended to inhomogeneous Dirichlet and Neumann boundary conditions.
Details on how inhomogeneous boundary conditions can be enforced are given in \cref{app:inhomogeneous-bcs}.

\section{Choice of Prior Covariance}\label{sec:choice-of-prior}
Thus far, the prior covariance matrix $\matSigma$ has not been specified.
The choice of $\matSigma$ is subject to two main requirements.
The first requirement is that $\matSigma$ needs to have a sparse representation.
Since $\matSigma$ is a $n \times n$ matrix, where $n$ is the number of degrees of freedom of the fine discretization, explicitly computing, storing and applying operations on the full matrix would quickly become prohibitively expensive.
As a result, the traditional approach of using a kernel to directly compute all entries of $\matSigma$ would be infeasible.
Instead, the prior is defined implicitly by assigning a sparse covariance matrix to the fine-scale force vector $\vecf$, which implicitly defines the covariance matrix of the solution vector $\matSigma$, but does not require us to explicitly compute it.
For certain kernel-based priors, an equivalent stochastic partial differential equation can be shown to exist, which allows for a similar sparse representation (see for example~\citet{roininen_whittle-matern_2014}).

The second requirement is that the choice of prior distribution needs to be appropriate for the partial differential equation at hand.
For instance, if the infinitely differentiable squared exponential prior were assumed on the solution field $u(x)$, this would imply $C^\infty$ continuity on the right-hand side field $f(x)$.
From a modeling point of view, this would be an undesirable assumption to make, since it is very restrictive concerning what forcing terms are permitted.
On the other hand, if the prior is not smooth enough, samples from the prior would exhibit unphysical discontinuities in $u(x)$ or its gradient fields.
In short, the prior needs to respect the smoothness of the partial differential equation to which it is applied.

\begin{figure*}
    \centering
    \includegraphics[width=0.6\textwidth]{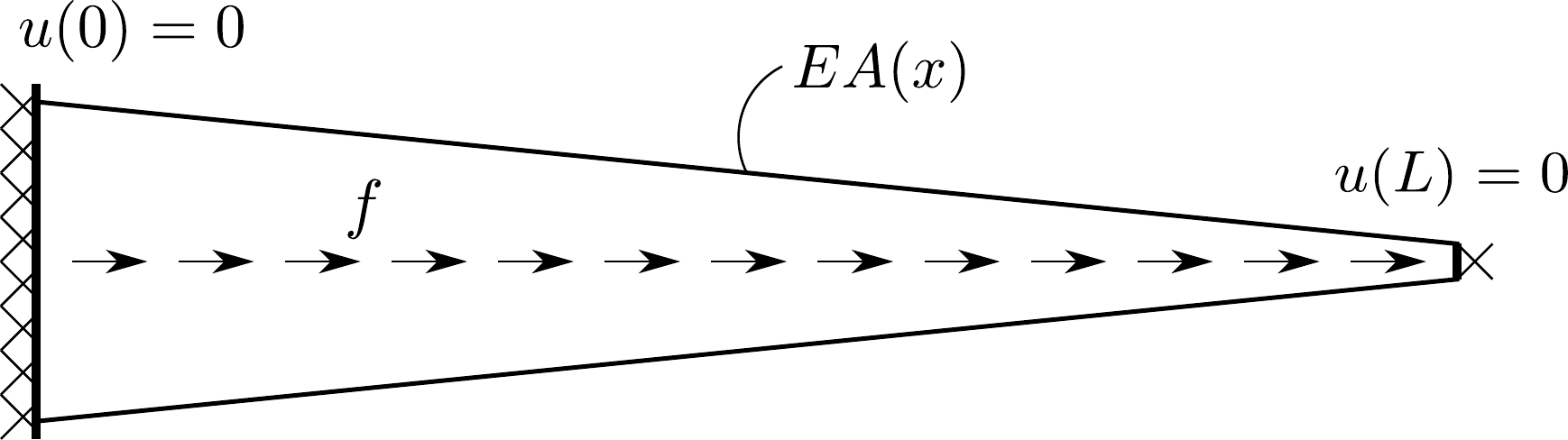}
    \caption{Schematic overview of the tapered bar test case}
    \label{fig:tapered-bar_overview}
\end{figure*}

In this section, a particular class of priors that meets both of these requirements is presented by means of two test cases.
The first test case, presented in \cref{fig:tapered-bar_overview}, concerns a one-dimensional mechanics problem described by the following ordinary differential equation with homogeneous boundary conditions:
\begin{equation}\label{eq:bar-strong-form}
\begin{aligned}
    -\frac{\diff}{\dx}\left(EA(x) \frac{\diff u}{\dx} \right) &= f(x) && \text{ in } \Omega = \left(0, 1\right)\\
    u(x) &= 0 && \text{ on } \partial \Omega = \{0,1\}
\end{aligned}
\end{equation}
Here, the distributed load $f(x) = 1$, Young's modulus $E = 1$ and the cross-sectional area $A(x) = 0.1 - 0.099 x$.
This setup describes a tapered bar with a constant load, where both the left and right end are clamped, as shown in \cref{fig:tapered-bar_overview}.
The fine-scale discretization consists of a uniform mesh with 64 elements ($n = 64$) and linear shape functions.
Three different levels of uniform coarse discretization are used: $m = 4$, $m = 16$ and $m = 64$.
Note that in all cases, since $n$ is a multiple of $m$, the shape functions are defined hierarchically in accordance with \cref{subsec:hierarchical-shape-functions}.

\begin{figure*}
    \centering
    \begin{subfigure}[b]{0.5\textwidth}
        \centering
        \includegraphics[width=\textwidth]{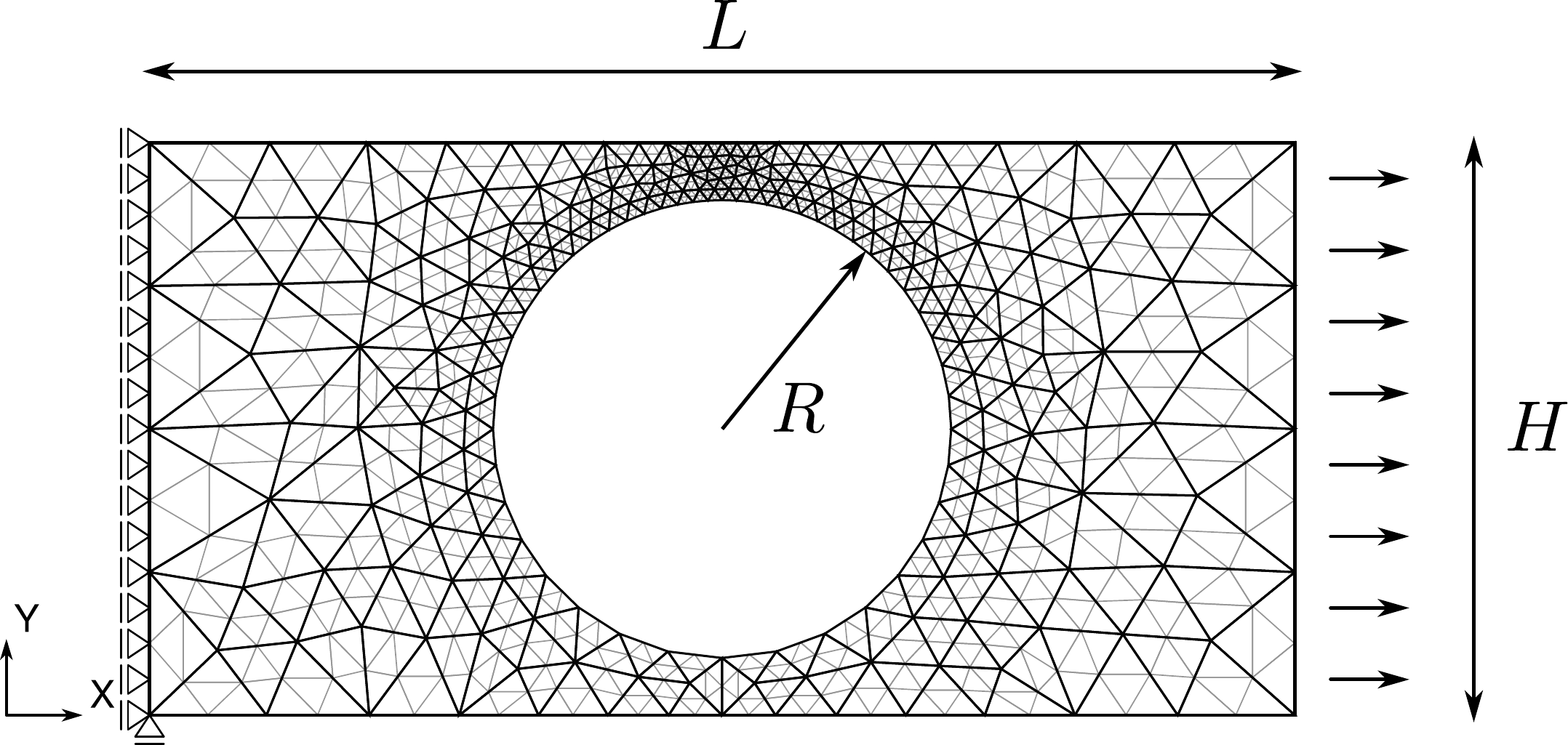}
        \subcaption{Schematic overview}
        \label{subfig:plate-with-hole_overview}
    \end{subfigure}\hfill
    \begin{subfigure}[b]{0.5\textwidth}
        \centering
        \includegraphics[width=0.9\textwidth]{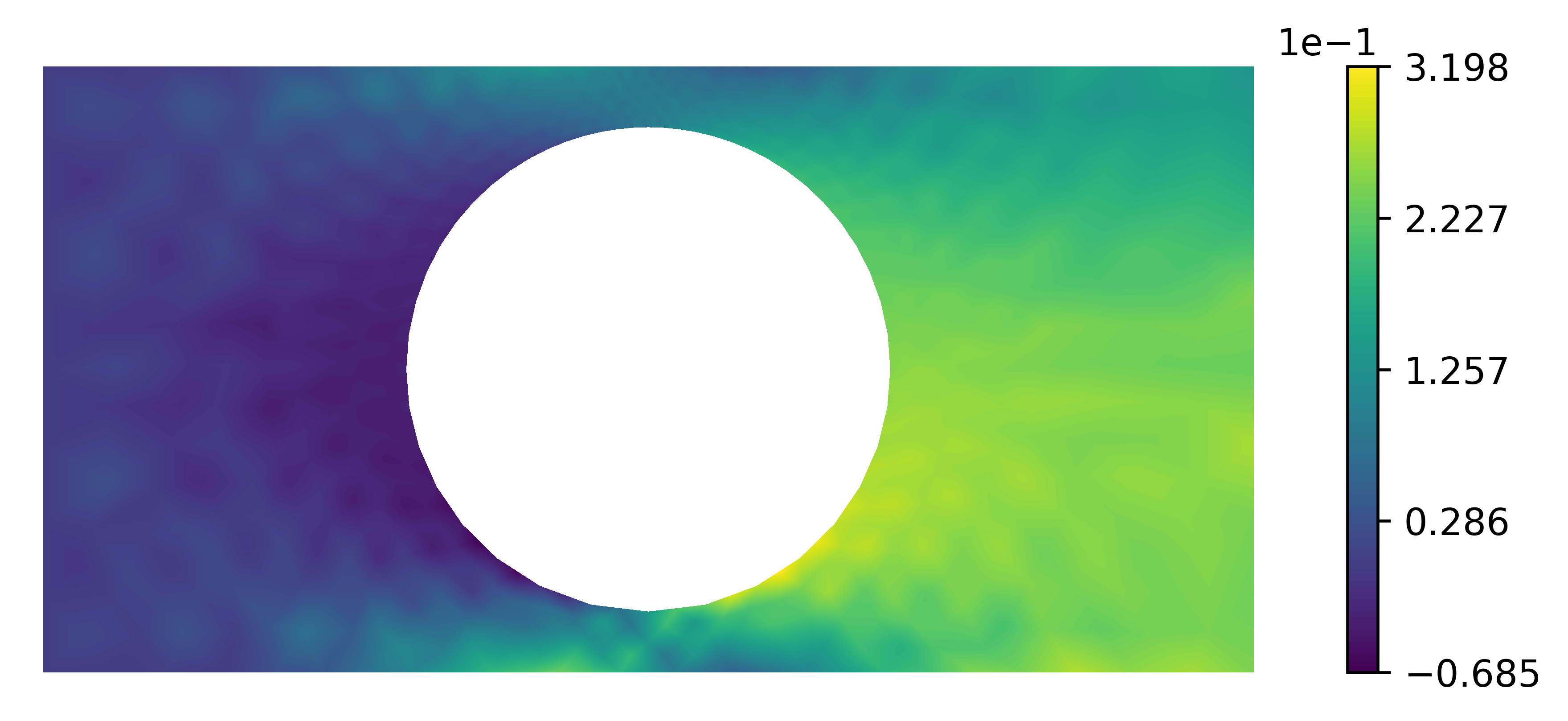}
        \subcaption{Discretization error $\vece$}
        \label{subfig:plate-with-hole_discretization-error}
    \end{subfigure}

    \begin{subfigure}[b]{0.5\textwidth}
        \centering
        \includegraphics[width=0.9\textwidth]{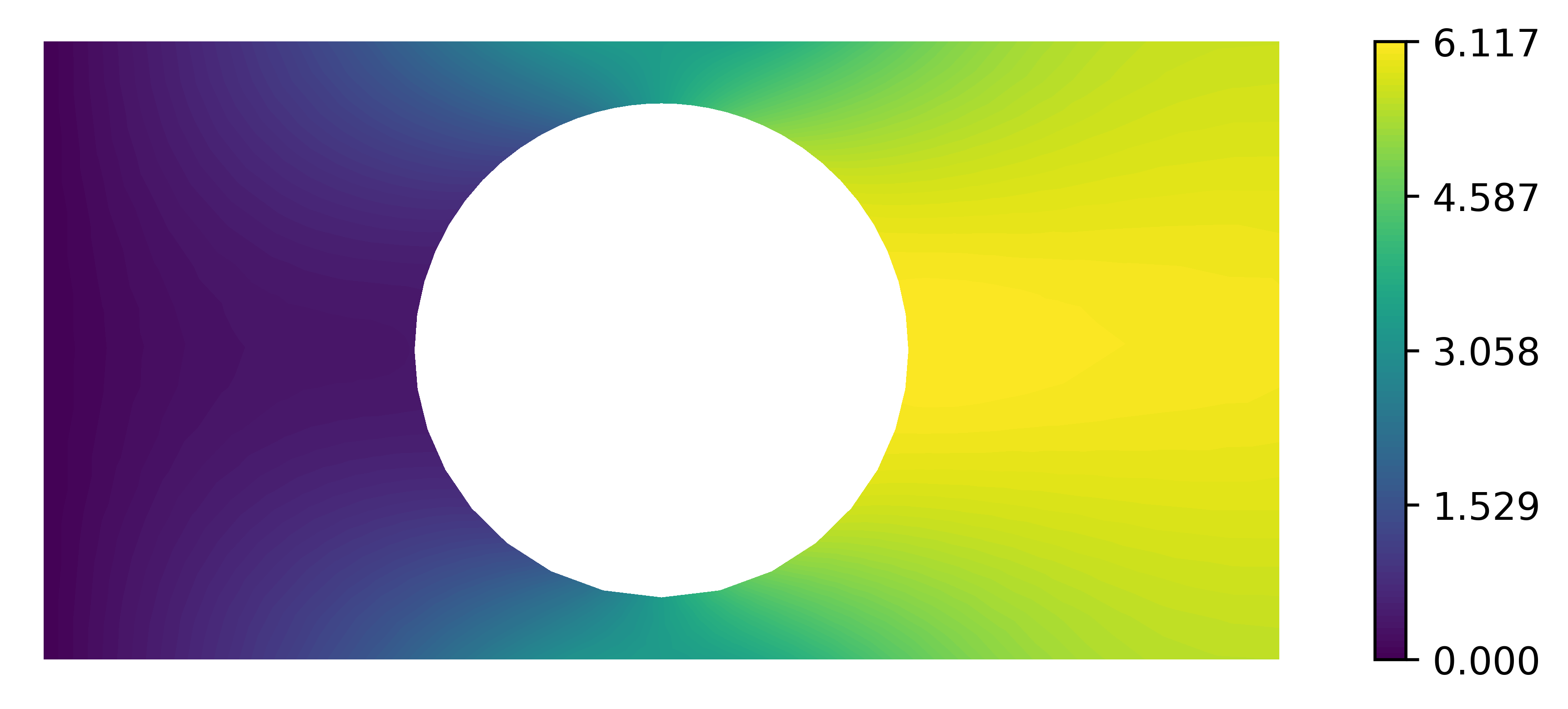}
        \subcaption{Fine-scale solution $\boldhat{\vecu}$}
        \label{subfig:plate-with-hole_fine-solution}
    \end{subfigure}\hfill
    \begin{subfigure}[b]{0.5\textwidth}
        \centering
        \includegraphics[width=0.9\textwidth]{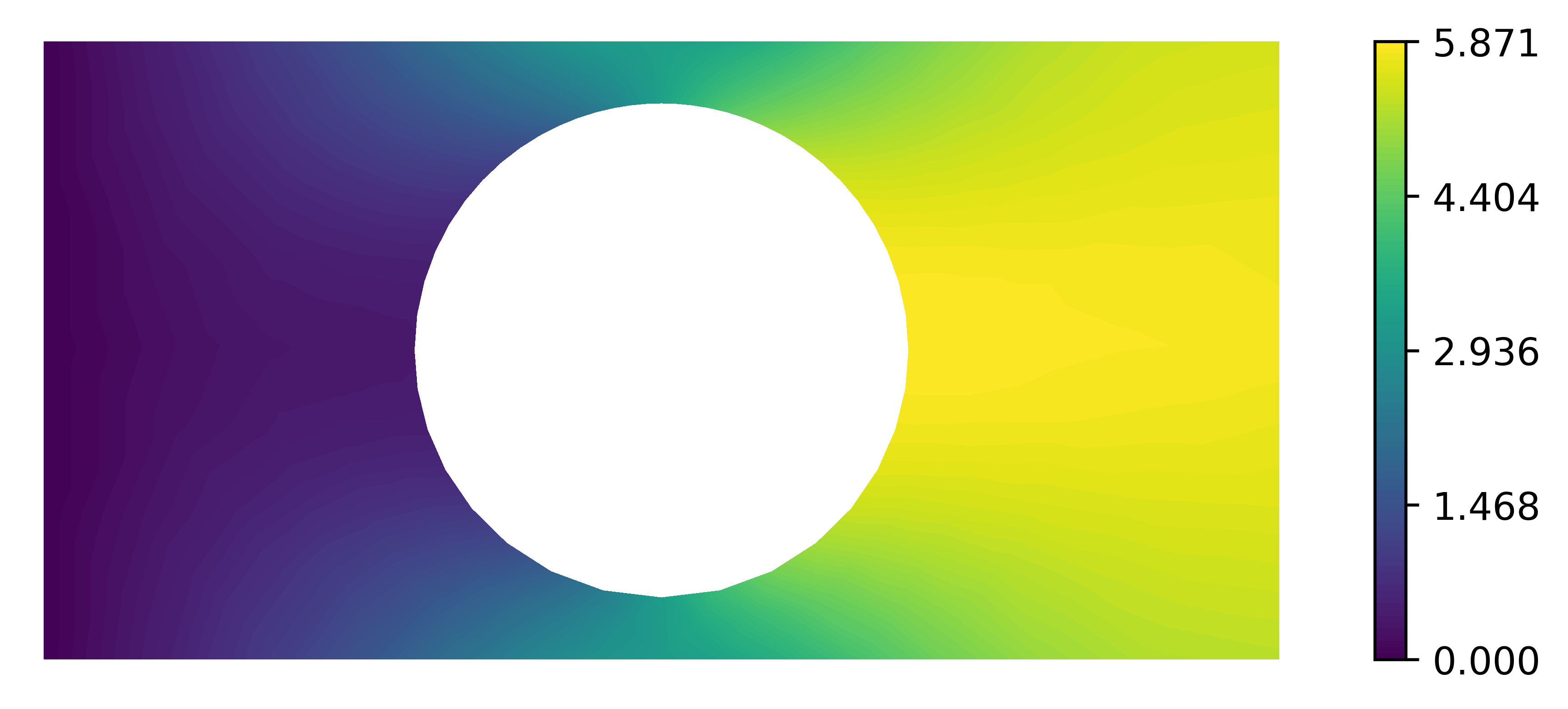}
        \subcaption{Coarse-scale solution $\vecuc$}
        \label{subfig:plate-with-hole_coarse-solution}
    \end{subfigure}
    \caption{Overview of the perforated plate bar test case}
    \label{fig:plate-with-hole_overview}
\end{figure*}

The second case is shown in \cref{fig:plate-with-hole_overview} and concerns a two-dimensional mechanics problem.
A plate ($L=4$, $H=2$) with a hole ($R = 0.8$) is clamped on its left edge and loaded by a constant horizontal body load $f_x = 1$
The plate has unit thickness, Young's modulus $E = 3$ and Poisson's ratio $\nu = 0.2$.
The problem is meshed non-uniformly, as shown in \cref{subfig:plate-with-hole_overview}.
For the coarse mesh, a characteristic length $h=0.5$ is used at the left and right edge, but around the hole a refinement is applied.
The refinements below and above the hole have a characteristic length of $h=0.2$ and $h=0.05$, respectively.
The fine mesh is generated by dividing each coarse element into 4 smaller triangular elements.
In \cref{subfig:plate-with-hole_discretization-error}, it can be seen how this difference in mesh density on different sides of the hole results in a larger discretization error below the hole than above it.
For reference, the fine-scale and coarse-scale solution are shown in \cref{subfig:plate-with-hole_fine-solution,subfig:plate-with-hole_coarse-solution}, respectively.

\subsection{A sparse right-hand side prior}\label{subsec:natural-prior}
Following the approach taken in~\citet{cockayne_probabilistic_2017}, rather than assuming a prior measure directly on the displacement field $u(\vecx)$, we assume a centered Gaussian process prior with covariance function $k_\text{f}(\vecx, \vecx')$ over the forcing term $f(\vecx)$:
\begin{equation}\label{eq:nat-prior-f-continuous}
\begin{aligned}
    f(\vecx) \sim \GP{0}{k_\text{f}(\vecx, \vecx')}
\end{aligned}
\end{equation}
This implicitly defines an equivalent prior on $u(\vecx)$:
\begin{equation}\label{eq:nat-prior-u-continuous}
\begin{aligned}
    u(\vecx) \sim \GP{0}{k_\text{nat}(\vecx, \vecx')}
\end{aligned}
\end{equation}
Here, the covariance function $k_\text{nat}$ can be expressed in terms of $k_\text{f}(\vecx, \vecx')$ and the Green's function $G(\vecx, \vecx')$ associated with the operator of the partial differential equation:
\begin{equation}\label{eq:natural-kernel}
\begin{aligned}
    k_\text{nat}(\vecx, \vecx') = \int_{\Omega} \int_{\Omega} G(\vecx, \vecz) G(\vecx', \vecz') k_\text{f}(\vecz, \vecz') \vecdz \vecdz'
\end{aligned}
\end{equation}
In~\citet{cockayne_probabilistic_2017}, this kernel is described as ``natural'' in the sense that the operator $-\Delta$ (see \cref{eq:strong-form}) uniquely maps from the Hilbert space associated with the forcing term covariance function $k_f(\vecx, \vecx')$ to the one associated with $k_\text{nat}(\vecx, \vecx')$.
Each sample from $u(\vecx)$ drawn from this natural kernel has an equivalent sample from $f(\vecx)$ and vice versa.
Unfortunately, since the Green's function is generally not available for a given partial differential equation, \citet{cockayne_probabilistic_2017} discards this natural kernel is then discarded in favor of a Matèrn or Wendland kernel with the appropriate level of smoothness.

However, because we avoid this problem by introducing the fine-scale discretization, there is no need here to step away from the natural prior approach.
Instead, it can be approximated by applying the fine-scale finite element discretization first, and only then finding the natural covariance matrix for the solution vector $\vecu$.
Given the prior distribution over $f(\vecx)$ in \cref{eq:nat-prior-f-continuous} and the definition of the force vector in \cref{eq:bg-stiffness-matrix}, it follows that:
\begin{equation}\label{eq:nat-prior-f-distribution}
\begin{aligned}
    \vecf \sim \normal{\vecnull}{\matSigmaf}
\end{aligned}
\end{equation}
where the force vector covariance matrix $\matSigmaf$ is given by:
\begin{equation}\label{eq:nat-prior-f-covariance}
\begin{aligned}
    \matSigmaf = \int_{\Omega} \int_{\Omega} k_\text{f}(\vecx, \vecx') \vecphi(\vecx) \vecphi(\vecx')^T \vecdx' \vecdx
\end{aligned}
\end{equation}
The resulting prior distribution over $\vecu$ then becomes:
\begin{equation}\label{eq:nat-prior-u-distribution}
\begin{aligned}
    \vecu \sim \normal{\vecnull}{\matK^{-1} \matSigmaf \matK^{-1}}
\end{aligned}
\end{equation}
Note the similarity to the natural kernel in \cref{eq:natural-kernel}, with $\matK^{-1}$ and $\matSigmaf$ taking a similar role as $G(\vecx, \vecx')$ and $k_\text{f}(\vecx - \vecx')$, respectively~\citep{peker_analyzing_2023}.
Also similarly, each sample from $\vecu$ has an equivalent sample from $\vecf$ and vice versa.
Conceptually, our choice of prior is the same as~\citet{cockayne_probabilistic_2017}, except that we are working in the finite-dimensional space of the discretized system, rather than the infinite-dimensional space of the original partial differential equation.
The advantage of working in the finite-dimensional space is that $\matK^{-1}$ is computable, and as a result the natural prior can still be used.

Given this choice of prior and using \cref{eq:phi-H-K-relation}, the posterior distribution of the displacement field is given by:
\begin{equation}\label{eq:nat-posterior-distribution}
\begin{aligned}
    \vecu | \vecg \sim \normal{\boldstar{\vecm}}{\boldstar{\matSigma}}
\end{aligned}
\end{equation}
with the following posterior mean $\boldstar{\vecm}$ and posterior covariance $\boldstar{\matSigma}$:
\begin{equation}\label{eq:nat-posterior-moments}
\begin{aligned}
    \boldstar{\vecm} &= \matK^{-1} \matSigmaf \matPhi \left(\matPhi^T \matSigmaf \matPhi\right)^{-1} \matPhi^T \vecf \\
    \boldstar{\matSigma} &= \matK^{-1} \left(\matI - \matSigmaf \matPhi \left(\matPhi^T \matSigmaf \matPhi \right)^{-1} \matPhi^T \right) \matSigmaf \matK^{-1}
\end{aligned}
\end{equation}

As mentioned earlier, since the posterior covariance matrix $\boldstar{\matSigma}$ is a full $n \times n$ matrix, explicitly computing it is computationally infeasible.
In~\cref{app:drawing-samples}, we explain how the prior and posterior distributions can be sampled without needing to assemble the full matrices $\matSigma^*$, $\matSigma$, $\matK^{-1}$ or $(\matPhi^T \matSigmaf \matPhi)^{-1}$.

\subsection{White noise prior}\label{subsec:white-noise-prior}
Within the natural prior framework, the main choice that remains is what right-hand side covariance function $k_\text{f}(\vecx, \vecx')$ to assume.
For now, we will follow the choice of \citet{cockayne_probabilistic_2017} to use the prior from \citet{owhadi_bayesian_2015}, and assume $k_\text{f}(\vecx, \vecx')$ to be a Dirac delta function $\delta(\vecx)$, scaled by a single hyperparameter $\alpha$:
\begin{equation}\label{eq:noise-prior-f-continuous}
\begin{aligned}
    k_\text{f}(\vecx, \vecx') = \alpha^2 \delta(\vecx - \vecx')
\end{aligned}
\end{equation}
This defines a white noise field over $f(\vecx)$ with a standard deviation that is equal to $\alpha$.
The covariance matrices $\matSigmaf$ and $\matSigma$ then follow directly from \cref{eq:nat-prior-f-covariance,eq:nat-prior-u-distribution}:
\begin{equation}\label{eq:noise-prior-covariance}
\begin{aligned}
    \matSigmaf &= \alpha^2 \matM \\
    \matSigma &= \alpha^2 \matK^{-1} \matM \matK^{-1}
\end{aligned}
\end{equation}
where $\matM$ is the fine-scale (square and symmetric) mass matrix, given by:
\begin{equation}\label{eq:bubnov-mass-matrix}
\begin{aligned}
    M_{ij} &= \int_{\Omega} \phi_i(\vecx) \phi_j(\vecx) \vecdx
\end{aligned}
\end{equation}
Note that under this choice of prior covariance, the sparsity requirement that was put on $\matSigma$ has been met.
The resulting posterior mean vector and covariance matrix are then given by:
\begin{equation}\label{eq:noise-posterior-moments}
\begin{aligned}
    \boldstar{\vecm} &= \matK^{-1} \matM \matPhi \left(\matPhi^T \matM \matPhi\right)^{-1} \matPhi^T \vecf \\
    \boldstar{\matSigma} &= \alpha^2 \matK^{-1} \left(\matI - \matM \matPhi \left(\matPhi^T \matM \matPhi \right)^{-1} \matPhi^T \right) \matM \matK^{-1}
\end{aligned}
\end{equation}
It can be seen that for this choice of prior, the hyperparameter $\alpha$ does not affect the posterior mean, and only serves as a scaling factor of the posterior covariance.
Given this hyperparameter-independence, we choose to simply set $\alpha = 1$ for the remainder of this work.
A small observation noise ($\sigma_e^2 = 10^{-12}$) is added to the term $\matPhi^T \matM \matPhi$ in \cref{eq:noise-posterior-moments}, to ensure that this matrix is invertible.

\begin{figure*}
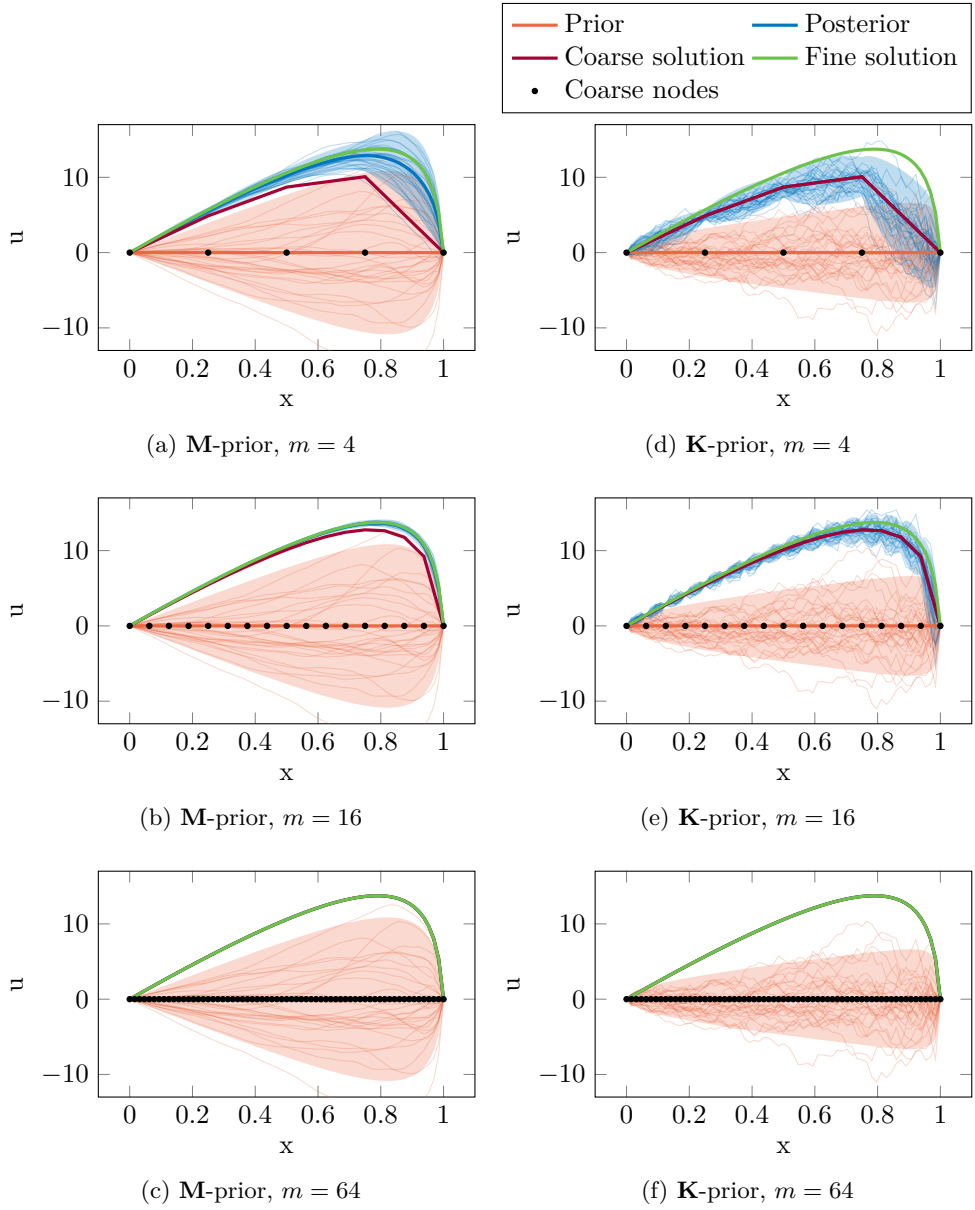

    \begin{subfigure}[b]{0.5\textwidth}
        \barplot{4}{0}{\textwidth}{0.7\textwidth}{M}
        \subcaption{$\matM$-prior, $m = 4$}
        \label{subfig:tapered-bar_prior-M_nc-4}
    \end{subfigure}\hfill
    \begin{subfigure}[b]{0.5\textwidth}
        \barplot{4}{1}{\textwidth}{0.7\textwidth}{K}
        \setcounter{subfigure}{3}
        \subcaption{$\matK$-prior, $m = 4$}
        \label{subfig:tapered-bar_prior-K_nc-4}
    \end{subfigure}

    \vspace{7pt}
    \begin{subfigure}[b]{0.5\textwidth}
        \barplot{16}{0}{\textwidth}{0.7\textwidth}{M}
        \setcounter{subfigure}{1}
        \subcaption{$\matM$-prior, $m = 16$}
        \label{subfig:tapered-bar_prior-M_nc-16}
    \end{subfigure}\hfill
    \begin{subfigure}[b]{0.5\textwidth}
        \barplot{16}{0}{\textwidth}{0.7\textwidth}{K}
        \setcounter{subfigure}{4}
        \subcaption{$\matK$-prior, $m = 16$}
        \label{subfig:tapered-bar_prior-K_nc-16}
    \end{subfigure}

    \vspace{7pt}
    \begin{subfigure}[b]{0.5\textwidth}
        \barplot{64}{0}{\textwidth}{0.7\textwidth}{M}
        \setcounter{subfigure}{2}
        \subcaption{$\matM$-prior, $m = 64$}
        \label{subfig:tapered-bar_prior-M_nc-64}
    \end{subfigure}\hfill
    \begin{subfigure}[b]{0.5\textwidth}
        \barplot{64}{0}{\textwidth}{0.7\textwidth}{K}
        \setcounter{subfigure}{5}
        \subcaption{$\matK$-prior, $m = 64$}
        \label{subfig:tapered-bar_prior-K_nc-64}
    \end{subfigure}
    \caption{Prior and posterior distributions of the 1D tapered bar problem. For comparison, the fine-scale and coarse-scale reference solutions have been included. From each distribution, 30 samples have been plotted. The shaded regions correspond to the 95\% credible intervals of the distributions.}
    \label{fig:tapered-bar_prior-M-K}
\end{figure*}

In \cref{subfig:tapered-bar_prior-M_nc-4}, the resulting prior and posterior distributions for the tapered bar problem are shown for the number of coarse elements $m$ equal to $4$.
Several pieces of information about the problem, in absence of knowledge of the right-hand side term, can be found encoded in the prior.
We can see how the enforcement of boundary conditions described in \cref{subsec:dirichlet-bcs} indeed results in a distribution whose samples respect the boundary conditions imposed at $x=0$ and $x=1$.
Furthermore, a larger prior standard deviation is found in the region where the bar is thinner, reflecting the fact that a small perturbation in the right-hand side in this region would have a more pronounced effect on the displacement field.
Considering the posterior distribution, we see that its mean falls between the coarse- and fine-scale reference solutions.
Lastly, it can be seen that the region where the posterior standard deviation is largest corresponds to the region where the discretization error is largest.

In \cref{subfig:tapered-bar_prior-M_nc-16,subfig:tapered-bar_prior-M_nc-64}, we have increased the number of degrees of freedom of the coarse mesh $m$ to $16$ and $64$ respectively, to study its effect on the posterior distribution.
As the coarse-scale solution approaches the fine-scale solution, the posterior mean approaches the fine-scale solution accordingly.
Additionally, the posterior standard deviation shrinks along with the discretization error until the coarse mesh density meets the fine one at $m = n = 64$.
At this point, only a small posterior standard deviation remains due to the small observation noise that was included in the model.

\begin{figure*}
    \begin{subfigure}{0.5\textwidth}
        \centering
        \includegraphics[width=0.8\textwidth]{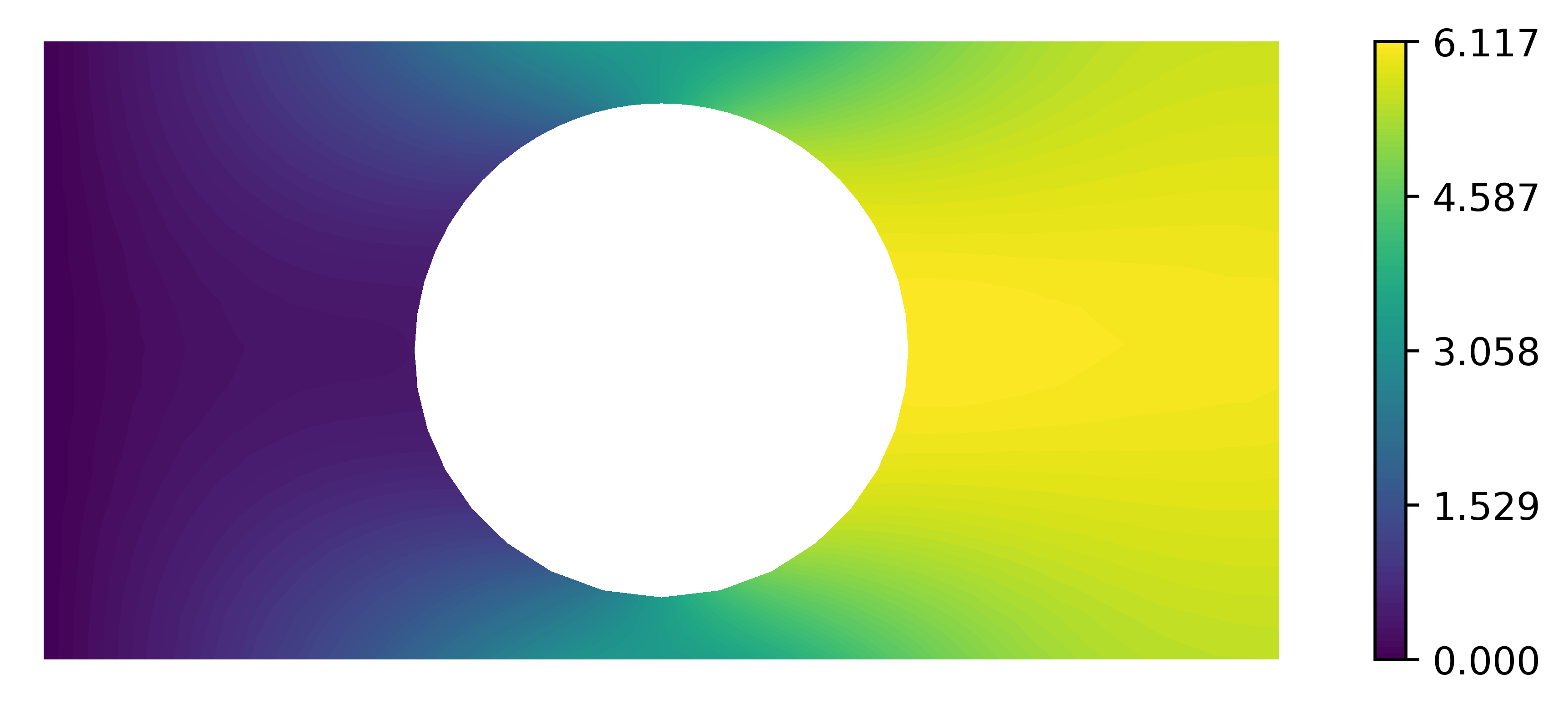}
        \subcaption{Posterior mean $\boldstar{\vecm}$}
        \label{subfig:plate-with-hole_prior-M_posterior-mean}
    \end{subfigure}\hfill
    \begin{subfigure}{0.5\textwidth}
        \centering
        \includegraphics[width=0.8\textwidth]{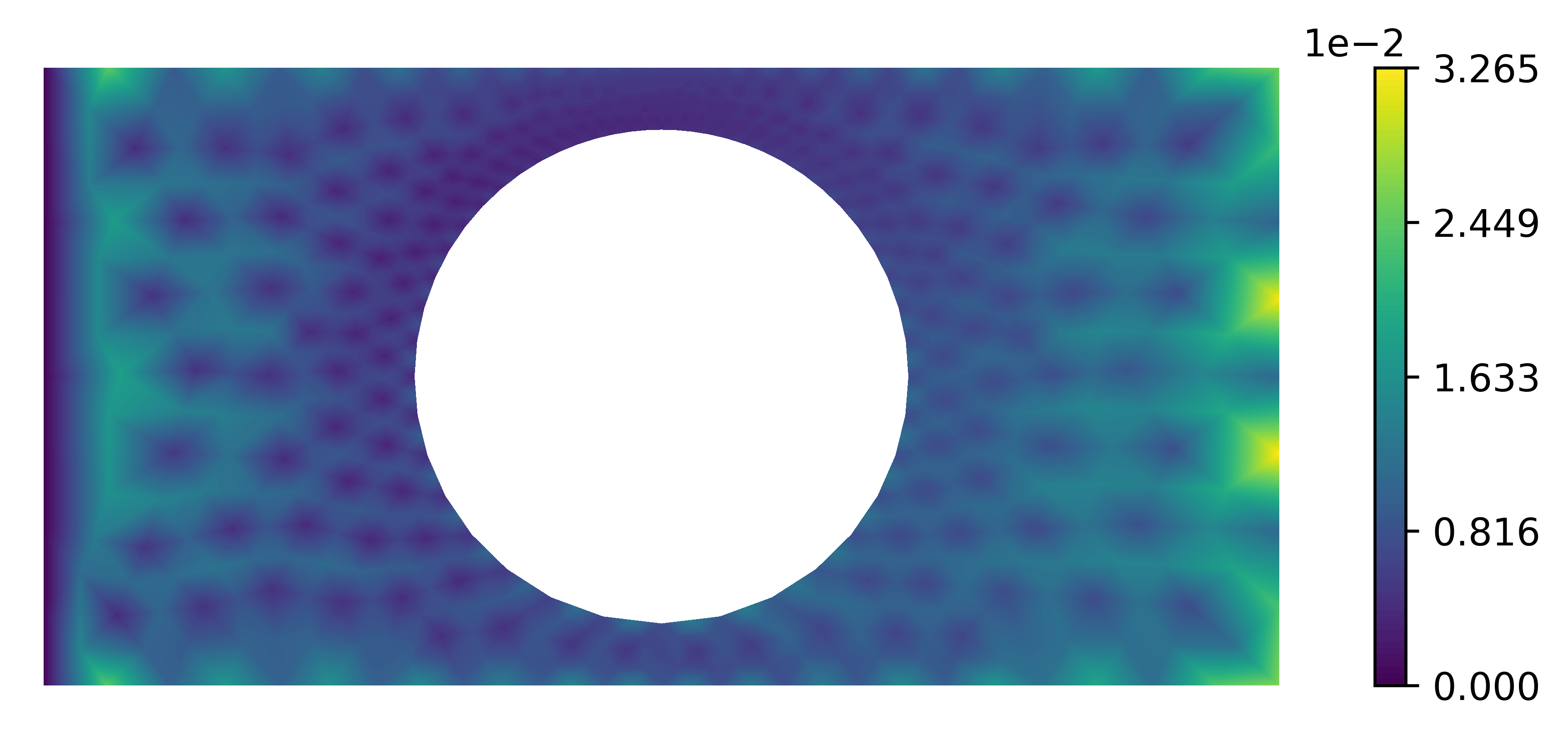}
        \subcaption{Posterior standard deviation $\boldstar{\vecsigma}$}
        \label{subfig:plate-with-hole_prior-M_posterior-std}
    \end{subfigure}
    \caption{Posterior moments of the perforated plate test case with $\matSigmaf = \matM$}
    \label{fig:plate-with-hole_prior-M}
\end{figure*}

In \cref{fig:plate-with-hole_prior-M}, the posterior moments are plotted when the same prior distribution is applied to the two-dimensional perforated plate problem.
We find that the results for this two-dimensional test case are quite different from those for the previous one-dimensional case.
It can be observed in \cref{subfig:plate-with-hole_prior-M_posterior-mean} that the posterior mean almost exactly matches the fine-scale solution shown in \cref{subfig:plate-with-hole_fine-solution}.
An explanation for this can be found by considering \cref{eq:nat-posterior-moments} and noting that the posterior mean $\boldstar{\vecm}$ is equivalent to the reference solution $\boldhat{\vecu}$, except that the force vector $\vecf$ has been replaced by $\boldhat{\vecf}$, a weighted projection of $\vecf$ onto the column space of $\matPhi$:
\begin{equation}\label{eq:f-hat}
\begin{aligned}
    \boldhat{\vecf} = \matP \vecf = \matSigmaf \matPhi \left(\matPhi^T \matSigmaf \matPhi \right)^{-1} \matPhi^T \vecf
\end{aligned}
\end{equation}
In other words, $\vecf$ is mapped to the coarse space, scaled, mapped back to the fine space and rescaled to obtain $\boldhat{\vecf}$.
The quality of this projection depends on the weights given by $\matSigmaf$, and there exists a sense in which the choice of $\matSigmaf = \matM$ is optimal:
it minimizes the projection error of the forcing term $f(\vecx)$ to the coarse space $\spaceWh$ in the $L^2$-norm (\citet{larson_finite_2013}, Theorem 1.1):
\begin{equation}\label{eq:optimal-l2-projection}
\argmin_{f^{\text{h}}(\vecx) \in \spaceWh} \| f(\vecx) - f^{\text{h}}(\vecx) \|_{L^2(\Omega)}^2 = \vecpsi(\vecx)^T \left(\matPhi^T \matM \matPhi \right)^{-1} \matPhi^T \vecf
\end{equation}
This optimality helps explain the close correspondence between the fine-scale solution and posterior mean given in \cref{subfig:plate-with-hole_fine-solution,subfig:plate-with-hole_prior-M_posterior-mean}, respectively.

Though this might appear to be a desirable property, for the purposes of modeling discretization error, it is actually detrimental.
To understand this, let us consider the following equality:
\begin{equation}\label{eq:contraction-measure}
\begin{aligned}
    \boldstar{\matSigma} \matSigma^{-1} \boldhat{\vecu} = \boldhat{\vecu} - \boldstar{\vecm}
\end{aligned}
\end{equation}
This expression can be easily verified by substituting the expressions for $\boldhat{\vecu}$, $\matSigma$, $\boldstar{\vecm}$ and $\boldstar{\matSigma}$ found in \cref{eq:nat-posterior-moments,eq:nat-prior-u-distribution,eq:reference-solution}.
The left-hand side of \cref{eq:contraction-measure} can be understood as quantifier of the amount of ``contraction'' of the prior distribution due to the observed data.
In the extreme case where there is no contraction of the covariance, we find that the posterior covariance matrix $\boldstar{\matSigma}$ is equal to the prior covariance matrix $\matSigma$, and consequently $\boldstar{\matSigma} \matSigma^{-1} = \matI$ and $\boldstar{\vecm} = \vecnull$.
At the other extreme, where the posterior covariance is given by $\boldstar{\matSigma} = \epsilon \matI$ and we let $\epsilon \to 0$, we find that the left-hand side approaches the null vector and as a result $\boldstar{\vecm} \to \boldhat{u}$.
As more observations are included, the posterior distribution moves from the former extreme to the latter.

It becomes clear that the posterior mean vector $\boldstar{\vecm}$ and posterior covariance matrix $\boldstar{\matSigma}$ are inextricably linked.
This property is not necessarily a problematic one.
In fact, from the typical probabilistic numerics point of view, where the solving procedure is interpreted as an inherently probabilistic process~\citep{hennig_probabilistic_2022}, the fact that the posterior covariance tends to zero as the posterior mean approaches the true solution is the desired kind of behavior.
However, if our goal is to have the discretization error reflected in the posterior covariance, then this connection to the posterior mean does pose a problem:
it is not possible to simultaneously obtain a posterior mean that approaches the true solution and a posterior covariance that is indicative of the coarse-scale discretization error.
And indeed, when comparing the posterior standard deviation $\vecsigma^*$ in \cref{subfig:plate-with-hole_prior-M_posterior-std} to the discretization error $\vece$ in \cref{subfig:plate-with-hole_discretization-error}, we see that the regions of largest discretization error are not reflected in the posterior standard deviation.

\subsection{Green's function prior}\label{subsec:greens-function-prior}
This crucial observation motivates us to reevaluate our initial choice of prior.
Given how \cref{eq:contraction-measure} relates the posterior covariance matrix $\boldstar{\matSigma}$ to the difference between the reference solution $\boldhat{\vecu}$ and the posterior mean vector $\boldstar{\vecm}$, it makes sense to choose a prior that will yield a posterior mean equal to the coarse-scale solution $\vecuc$.
Additionally, from a discretization error modeling point of view, it is more sensible to have a posterior mean that is equal to the coarse-scale solution $\vecuc$ than to have a posterior mean that improves on it.
After all, the aim from the outset has been to interpret the finite element discretization error as a source of uncertainty surrounding the coarse-scale finite element solve.

In~\citet{pfortner_physics-informed_2023}, a method is presented to construct a prior whose posterior mean matches exactly the coarse-scale finite element solution $\vecuc$ from an initial prior an arbitrary mean function $m(\vecx)$ and covariance function $k(\vecx, \vecx')$.
However, we will opt instead for the method presented in~\citet{bilionis_probabilistic_2016}, which is to set the prior covariance function equal to the Green's function $G(\vecx, \vecx')$ of the partial differential equation at hand.
For Poisson's equation, this choice of prior yields the following right-hand side covariance function $k_\text{f}(\vecx, \vecx')$:
\begin{equation}\label{eq:greens-prior-f-continuous}
\begin{aligned}
    k_\text{f}(\vecx, \vecx') = - \Delta \delta(\vecx - \vecx')
\end{aligned}
\end{equation}
Substitution of this expression into \cref{eq:nat-prior-f-covariance}, applying integration by parts and subsequent substitution into \cref{eq:nat-prior-u-distribution} yields the following expressions for $\matSigmaf$ and $\matSigma$:
\begin{equation}\label{eq:greens-prior-covariance}
\begin{aligned}
    \matSigmaf &= \matK \\
    \matSigma &= \matK^{-1}
\end{aligned}
\end{equation}
Intuitively, we again find $\matK^{-1}$ as the finite-dimensional counterpart of $G(\vecx, \vecx')$.
The advantages of introducing the fine-scale discretization as a stand-in for the infinite-dimensional partial differential equation once again become apparent:
the fact that Green's function is generally unavailable does not pose a problem anymore.
Furthermore, the objection raised in~\citet{bilionis_probabilistic_2016} that for Poisson's equation in two or three dimensions, the Green's function $G(\vecx, \vecx')$ is infinite at $\vecx = \vecx'$ and can therefore not be a useful indicator of model uncertainty does not apply in our case:
for any valid finite element discretization the finite-dimensional inverse stiffness matrix $\matK^{-1}$ only has finite-valued entries.
In the phrasing of~\citet{alberts_physics-informed_2023}, the introduction of the fine-scale discretization offers a way to truncate the integration over functions at the smallest scales.

This choice of prior in turn results in the following posterior mean vector and covariance matrix:
\begin{equation}\label{eq:greens-posterior-moments}
\begin{aligned}
    \boldstar{\vecm} &= \matPhi \left(\matPhi^T \matK \matPhi\right)^{-1} \matPhi^T \vecf \\
    \boldstar{\matSigma} &= \matK^{-1} - \matPhi \left(\matPhi^T \matK \matPhi \right)^{-1} \matPhi^T
\end{aligned}
\end{equation}
Note that according to \cref{eq:phi-Kc-K-relation,eq:phi-g-f-relation}, $\matPhi^T \matK \matPhi$ and $\matPhi^T \vecf$ are equal to the coarse stiffness matrix $\matKc$ and coarse force vector $\vecg$, respectively.
As a result, we find that indeed the posterior mean vector $\boldstar{\vecm}$ is exactly equal to the solution of the coarse system $\vecuc$, projected to the fine space.
Returning now to \cref{eq:contraction-measure}, we find that for this choice of prior, this expression simplifies to a surprisingly simple relationship between the posterior covariance matrix $\boldstar{\matSigma}$ and the discretization error $\vece$ as defined in \cref{eq:discretization-error-definition}:
\begin{equation}\label{eq:greens-error-approximation}
\begin{aligned}
    \boldstar{\matSigma} \vecf = \vece
\end{aligned}
\end{equation}
Since this relation holds for any fine-scale force vector $\vecf$, and $\boldstar{\matSigma}$ is independent of $\vecf$, this posterior covariance matrix $\boldstar{\matSigma}$ can be used to determine the discretization error $\vece$ for an arbitrary forcing term.
In this sense, $\boldstar{\matSigma}$ can be said to fully encode the discretization error associated with the geometry and discretization of the problem at hand.

In \cref{subfig:tapered-bar_prior-K_nc-4,subfig:tapered-bar_prior-K_nc-16,subfig:tapered-bar_prior-K_nc-64}, the prior and posterior distributions that follow when applying this prior to the tapered bar problem are shown.
Again, the number of coarse elements $m$ is equal to $4$, $16$ and $64$, respectively.
As expected, the posterior mean $\boldstar{\vecm}$ can be seen to equal the coarse-scale solution $\vecuc$ in these figures.
Similar to \cref{subfig:tapered-bar_prior-M_nc-4,subfig:tapered-bar_prior-M_nc-16,subfig:tapered-bar_prior-M_nc-64}, the largest posterior standard deviation $\boldstar{\vecsigma}$ is found in the region where the bar is thinnest.
However, for this prior there is a notable reduction of the posterior standard deviation around the coarse-scale nodes.
This reduction in the standard deviation is reflective of the fact that at these nodes, the coarse solution is more accurate than in the regions between the coarse-scale nodes, where the solution is interpolated via the coarse-scale shape functions.

Another notable difference between the two priors in \cref{fig:tapered-bar_prior-M-K} is the smoothness of the samples.
We see in \cref{subfig:tapered-bar_prior-M_nc-4,subfig:tapered-bar_prior-M_nc-16,subfig:tapered-bar_prior-M_nc-64} that the samples from the white noise prior presented in \cref{subsec:white-noise-prior} have a visible smoothness to them.
In contrast, the samples from the Green's function prior shown in \cref{subfig:tapered-bar_prior-K_nc-4,subfig:tapered-bar_prior-K_nc-16,subfig:tapered-bar_prior-K_nc-64} appear jagged and rough.
In fact, for this one-dimensional Poisson problem, each sample $\tilde{u}(x)$ drawn from the Green's function prior $k(x, x') = G(x, x')$ can be shown to be continuous, but nowhere differentiable.
This is the result of the fact that at $x = x'$, the Green's function is continuous (i.e.\ $\lim_{\delta \to 0} G(x-\delta, x) = \lim_{\delta \to 0} G(x+\delta, x)$), but at that same point its derivative is discontinuous (i.e.\ $\lim_{\delta \to 0} G'(x - \delta, x) \neq \lim_{\delta \to 0} G'(x + \delta, x)$)~\citep{bayin_mathematical_2006}.
The samples of a Gaussian process are mean-square continuous if $k(\vecx, \vecx')$ is continuous at $\vecx = \vecx'$ and are $k$ times mean-square differentiable if $k(\vecx, \vecx')$ is $2 k$ times differentiable at $\vecx = \vecx'$~\citep{rasmussen_gaussian_2005}.
From the fact that the Green's function is not differentiable at $x = x'$, it thus follows that the samples drawn from this process are everywhere continuous but nowhere differentiable.
Note that this only applies to the infinite-dimensional solution space $\spaceV$.
The finite-dimensional space $\spaceVh$ spanned by the fine-scale shape functions $\vecphi(x)$ is still weakly once-differentiable for both priors.

\begin{figure*}
    \begin{subfigure}{0.5\textwidth}
        \centering
        \includegraphics[width=0.8\textwidth]{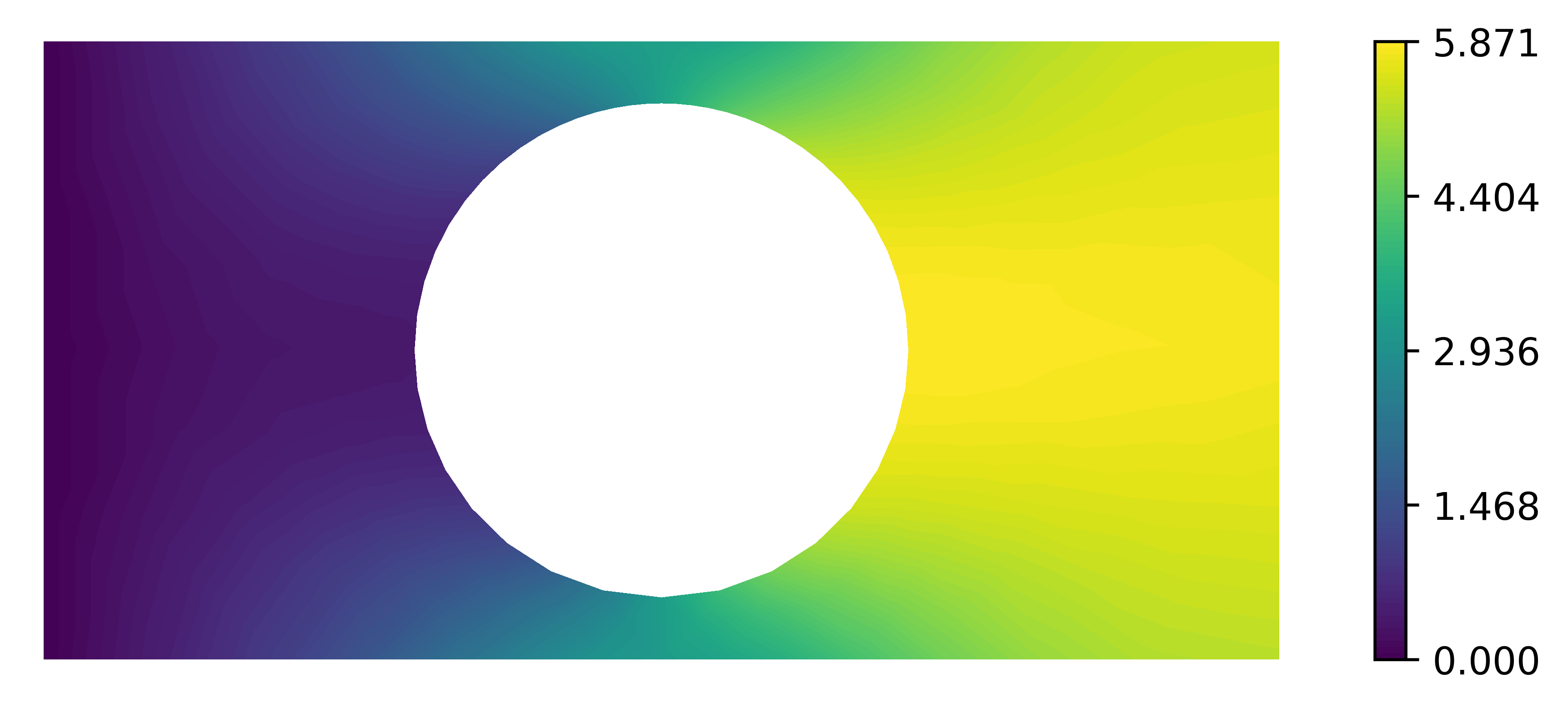}
        \subcaption{Posterior mean $\boldstar{\vecm}$}
        \label{subfig:plate-with-hole_prior-K_posterior-mean}
    \end{subfigure}\hfill
    \begin{subfigure}{0.5\textwidth}
        \centering
        \includegraphics[width=0.8\textwidth]{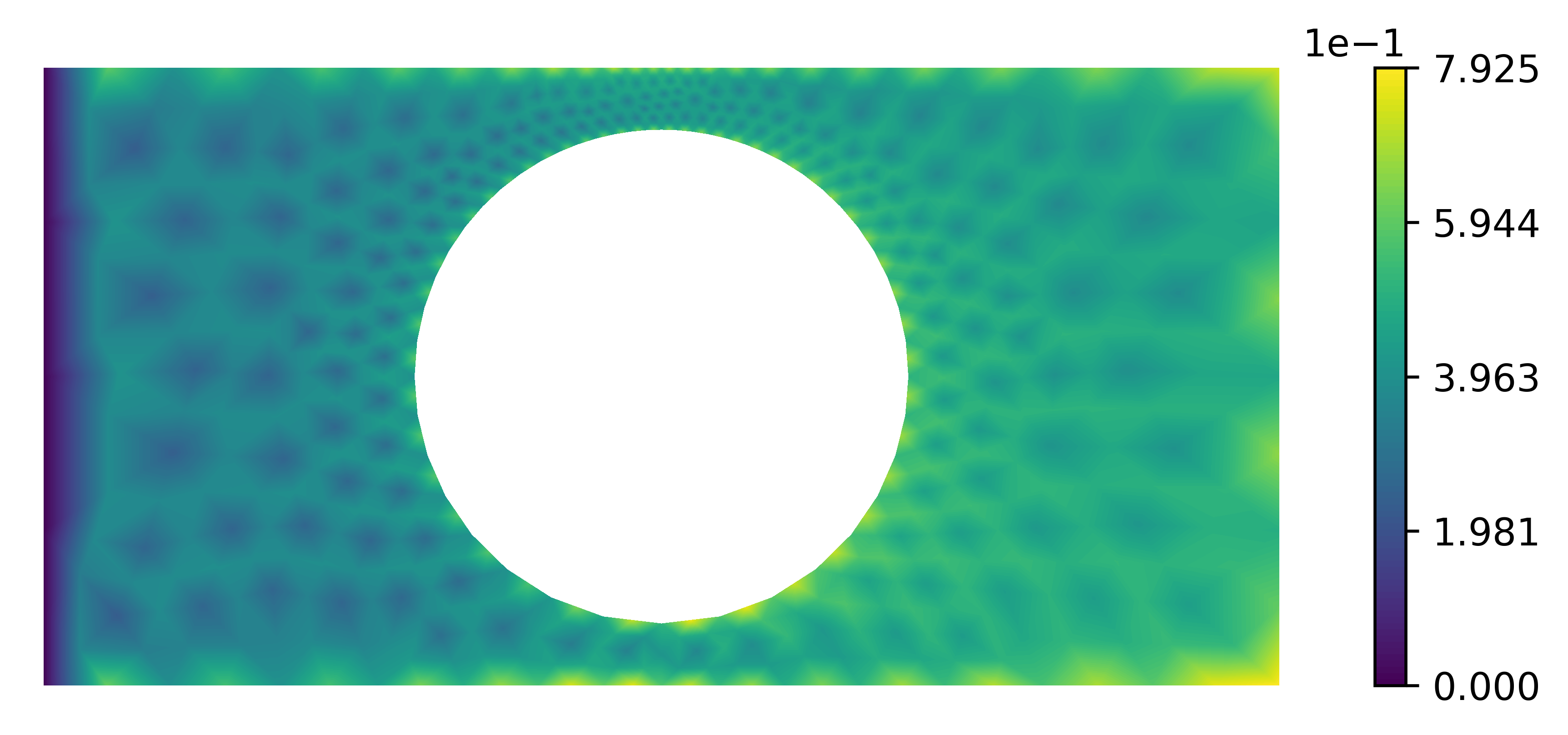}
        \subcaption{Posterior standard deviation $\boldstar{\vecsigma}$}
        \label{subfig:plate-with-hole_prior-K_posterior-std}
    \end{subfigure}

    \begin{subfigure}{0.5\textwidth}
        \centering
        \includegraphics[width=0.8\textwidth]{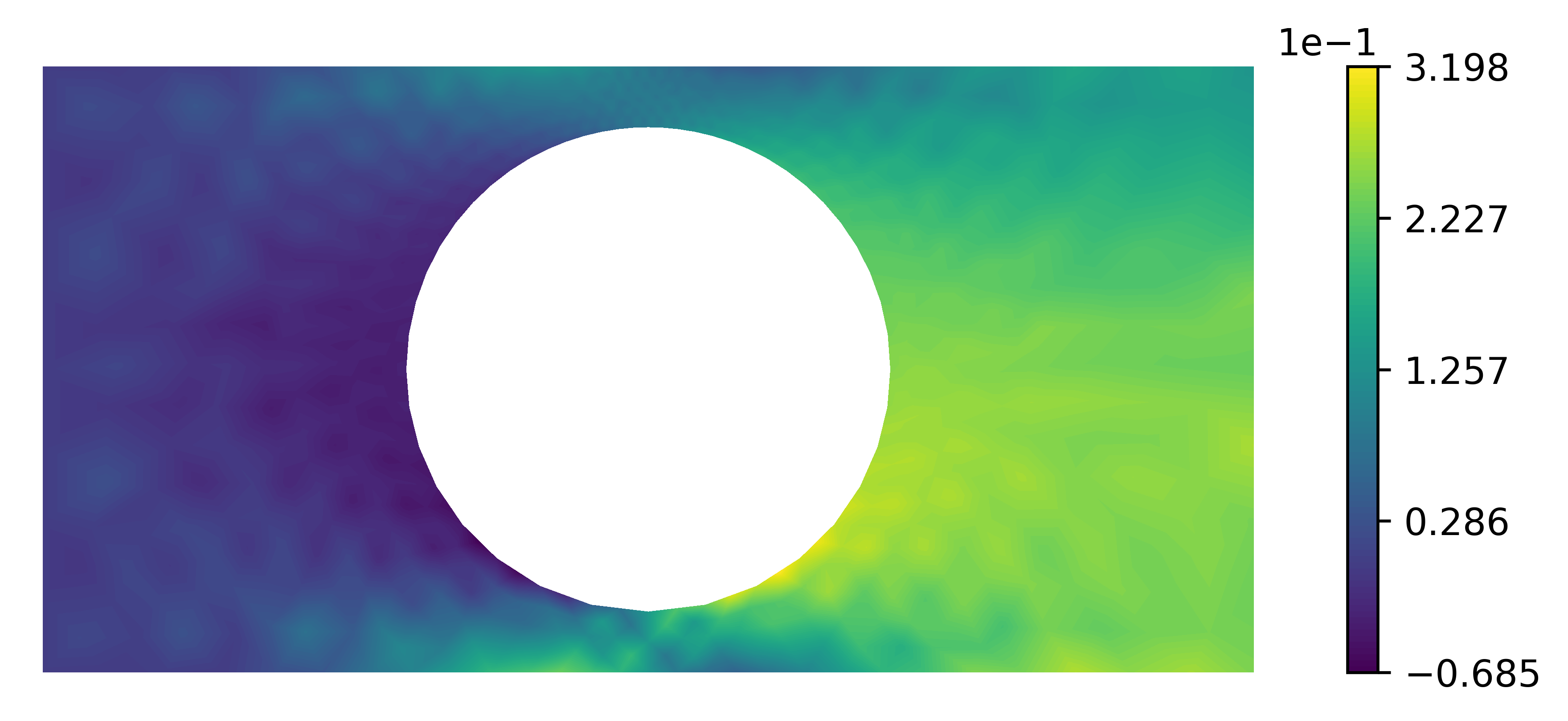}
        \subcaption{Error recovery $\boldstar{\matSigma} \vecf$}
        \label{subfig:plate-with-hole_prior-K_error-recovery}
    \end{subfigure}\hfill
    \begin{subfigure}{0.5\textwidth}
        \centering
        \includegraphics[width=0.8\textwidth]{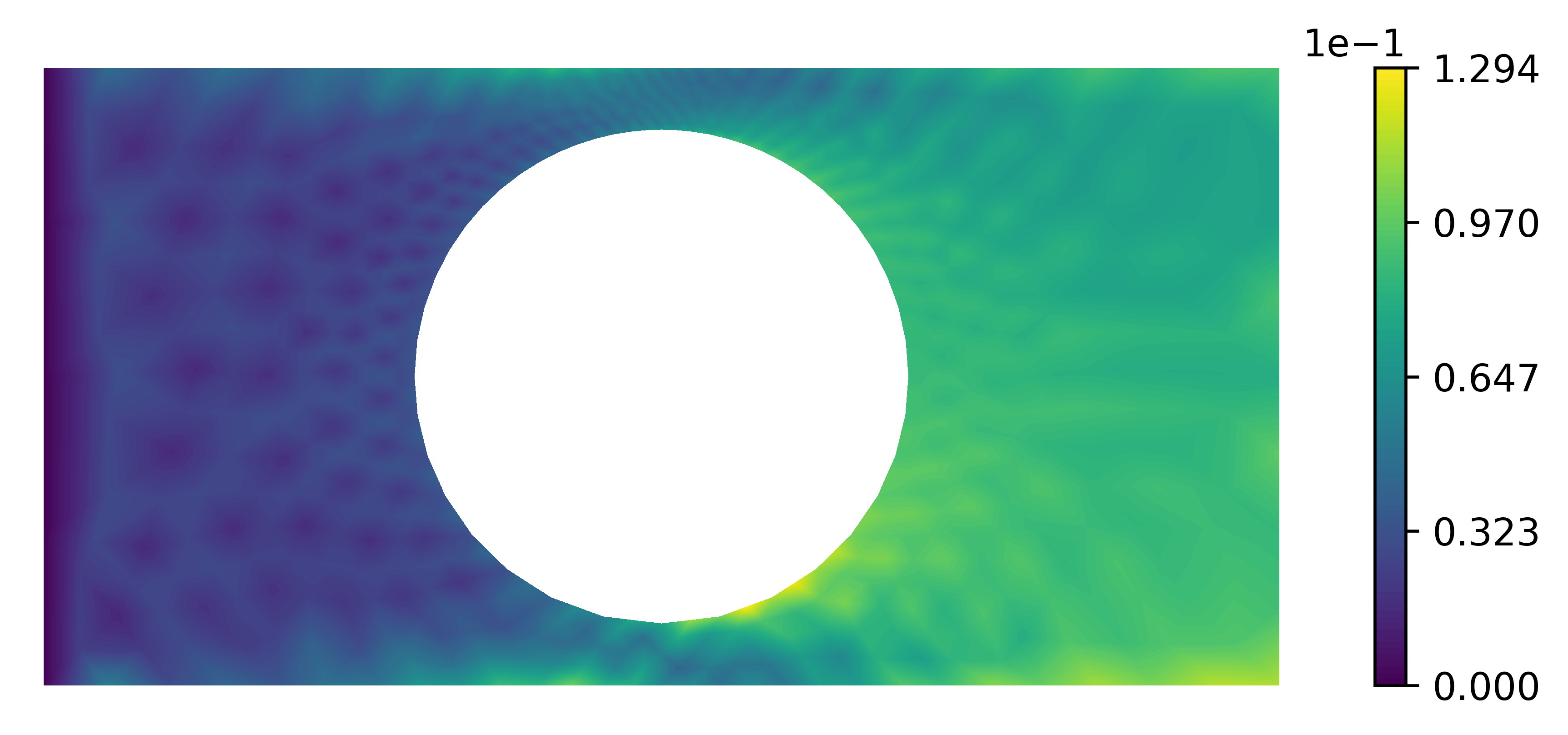}
        \subcaption{Rescaled posterior standard deviation $\boldstar{\boldhat{\vecsigma}}$}
        \label{subfig:plate-with-hole_prior-K_posterior-std-rescaled}
    \end{subfigure}
    \caption{Posterior moments of the perforated plate test case with $\matSigmaf = \matK$}
    \label{fig:plate-with-hole_prior-K}
\end{figure*}

We now turn to the perforated plate example, for which the results are shown in \cref{fig:plate-with-hole_prior-K}.
The posterior mean $\boldstar{\vecm}$ in \cref{subfig:plate-with-hole_prior-K_posterior-mean} can be seen to exactly match the coarse-scale finite element solution $\vecuc$ in \cref{subfig:plate-with-hole_coarse-solution} for this problem as well.
Unfortunately, the posterior covariance $\boldstar{\vecsigma}$ shown in \cref{subfig:plate-with-hole_prior-K_posterior-std} appears again to bear little resemblance to the discretization error $\vece$ from \cref{subfig:plate-with-hole_discretization-error}.
This might seem surprising, given the direct relationship between posterior covariance $\boldstar{\matSigma}$ and discretization error $\vece$ given in \cref{eq:greens-error-approximation}.
Indeed, we can multiply the posterior covariance $\boldstar{\matSigma}$ by the fine-scale force vector $\vecf$ to recover the discretization error exactly (see \cref{subfig:plate-with-hole_prior-K_error-recovery}), but this does not translate to a posterior standard deviation $\boldstar{\vecsigma}$ that can be interpreted directly.
This is a consequence of the fact that the posterior covariance $\boldstar{\matSigma}$ depends only on the material stiffness (via $\matK$) and node locations (via $\matPhi$), but not on the magnitude of the force vector $\vecf$ at those locations.
One benefit that results from this independence is that given the posterior covariance matrix $\boldstar{\matSigma}$ from one load case, it is possible to compute the discretization error for any other load case virtually for free\footnote{All that is needed is a matrix-vector multiplication of the posterior covariance matrix $\boldstar{\matSigma}$ and the fine-scale force vector $\vecf$ of the new load case}.
However, the drawback of this independence is that, since the discretization error $\vece$ does depend on the load applied to the structure, a load-independent posterior standard deviation $\boldstar{\vecsigma}$ cannot adequately represent the discretization error for any specific load case.
Paradoxically, because the posterior covariance matrix $\boldstar{\matSigma}$ encodes the discretization error $\vece$ for all load cases simultaneously, it fails to represent the discretization error for any one load case in particular.
This paradox is not unique to our Bayesian formulation of the finite element method, and arises in many Gaussian process--based probabilistic solver of differential equations, including meshfree probabilistic solvers~\citep{bilionis_probabilistic_2016, cockayne_probabilistic_2017} and probabilistic methods of weighted residuals~\citep{pfortner_physics-informed_2023}.
In all these cases, the error between the posterior mean function $m^*(\vecx)$ and exact solution $u(\vecx)$ is dependent on the right-hand side term, but the posterior covariance function $k^*(\vecx, \vecx')$ meant to represent this error is not.

\subsection{Incorporating force term information}\label{subsec:pca}
This raises the question whether it is possible to break this independence of the posterior covariance matrix $\boldstar{\matSigma}$ and the fine-scale force vector $\vecf$.
Doing so appears to be necessary to capture the load-dependent discretization error $\vece$ in the posterior standard deviation $\boldstar{\vecsigma}$.
Returning to \cref{eq:greens-error-approximation}, we can understand the multiplication of $\boldstar{\matSigma}$ by $\vecf$ through the eigendecomposition of $\boldstar{\matSigma}$:
\begin{equation}\label{eq:posterior-eigendecomposition}
\begin{aligned}
    \boldstar{\matSigma} = \matQ \matLambda \matQ^{-1}
\end{aligned}
\end{equation}
Here the columns of $\matQ$ are the eigenvectors of $\boldstar{\matSigma}$ and $\matLambda$ is a diagonal matrix whose entries are its eigenvalues in descending order.
Since $\boldstar{\matSigma}$ is real positive definite, its eigenvalues are all positive real numbers, and $\matQ$ is an orthogonal matrix, which implies that $\matQ^{-1} = \matQ^T$.

The decomposition in \cref{eq:posterior-eigendecomposition} allows for a straightforward interpretation of the multiplication of $\boldstar{\matSigma}$ by $\vecf$.
First, $\matQ^{-1}$ performs a change of basis $\boldtilde{\vecf} = \matQ^{-1} \vecf$, expressing $\vecf$ in terms of the basis spanned by the eigenvectors instead of the standard basis.
In this basis, $\boldtilde{\vecf}$ is rescaled by the eigenvalues $\matLambda$ to obtain the discretization error $\boldtilde{\vece}$ expressed in terms of the eigenbasis.
Finally, $\matQ$ performs a change of basis on $\boldtilde{\vece}$ back to the standard basis $\vece = \matQ \boldtilde{\vece}$.
Since $\matLambda$ is a diagonal matrix, the operation $\boldtilde{\vece} = \matLambda \boldtilde{\vecf}$ comes down to a simple element-wise multiplication:
\begin{equation}\label{eq:eigenbasis-error-approximation}
\begin{aligned}
    \tilde{e}_i = \lambda_i \tilde{f}_i
\end{aligned}
\end{equation}
Rather than interpreting \cref{eq:eigenbasis-error-approximation} as a rescaling of each element of the force vector $\tilde{f}_i$ by its corresponding eigenvalue $\lambda_i$, one could argue equally well that it is the eigenvalue $\lambda_i$ that is rescaled by $\tilde{f}_i$ instead.
In order to break the independence of the posterior covariance matrix $\boldstar{\matSigma}$ and force vector $\vecf$, we replace the original eigenvalues $\lambda_i$ with ones that are rescaled by $\tilde{f}_i$
Thus, $\matLambda$ is replaced by a diagonal matrix $\matE$, whose diagonal entries are given by $|\tilde{e}_i|$, yielding a new covariance matrix $\boldstar{\boldhat{\matSigma}}$:
\begin{equation}\label{eq:rescaled-posterior-covariance}
\begin{aligned}
    \boldstar{\boldhat{\matSigma}} = \matQ \matE \matQ^{-1}
\end{aligned}
\end{equation}
Since all entries of $\matE$ are nonnegative, this rescaled covariance matrix $\boldstar{\boldhat{\matSigma}}$ is positive semi-definite, and thus a valid covariance matrix.
In \cref{subfig:plate-with-hole_prior-K_posterior-std-rescaled}, the standard deviation $\boldstar{\boldhat{\vecsigma}}$ of this rescaled covariance matrix is shown.
Comparing to the discretization error $\vece$ in \cref{subfig:plate-with-hole_discretization-error}, we see a clear similarity between these two fields.
At last, we appear to have arrived at a distribution with a covariance matrix that can meaningfully capture the discretization error.

One shortcoming of this ad hoc approach to incorporating forcing term information in our posterior distribution, is that it is a deviation from the Bayesian paradigm used thus far, since there is no guarantee that there exists an equivalent prior distribution that would yield this rescaled posterior covariance matrix.
Additionally, if there does exist an equivalent prior, it is unclear what posterior mean this equivalent prior would produce.
Our motivation for presenting this approach nonetheless is to demonstrate not only that it is impossible to obtain an interpretable posterior standard deviation $\boldstar{\vecsigma}$ if the posterior covariance matrix $\boldstar{\matSigma}$ is independent of the forcing term $\vecf$, but also that it is possible to obtain an interpretable standard deviation by incorporating forcing term information.

\section{Conclusions}\label{sec:conclusions}
In this work, we presented a Bayesian approach to the modeling of finite element discretization error.
A Gaussian process prior is assumed over the solution space, which is conditioned on the force vector from a finite element discretization.
To avoid the computation of intractable integrals, a second, finer mesh is introduced, which is assumed to be sufficiently fine to represent the true solution.
The two meshes are constructed in a hierarchical manner, such that the coarse-scale shape functions can be fully expressed in terms of fine-scale shape functions.
The Gaussian process prior on the solution space yields a normal distribution prior on the fine-scale solution vector.
For linear partial differential equations, conditioning this prior on the coarse-scale force vector produces a normally distributed posterior on the solution vector.

Two different prior covariance functions have been investigated: a white noise prior covariance on the forcing term, and a Green's function prior covariance on the solution term.
The white noise prior covariance is shown to produce a posterior mean vector that is close to the fine-scale reference solution.
However, an undesirable consequence of this property is that the corresponding posterior covariance matrix becomes less informative of the discretization error between the coarse-scale and fine-scale solutions.
This is a result of the fact that the posterior covariance matrix is inherently connected to the difference between the posterior mean and fine-scale solution.
As a consequence, in order to obtain a posterior covariance that is informative of the difference between the coarse-scale and fine-scale solutions, one must choose a prior that will produce a posterior mean equal to the coarse-scale solution.
This results in a different choice of prior than most probabilistic numerical methods, where the aim is to converge to the exact solution.

This leads us to the Green's function prior covariance, which is shown to produce exactly the coarse-scale solution as its posterior mean.
Additionally, the discretization error can be recovered exactly from the posterior covariance matrix by multiplying it by the fine-scale force vector.
Because the posterior covariance matrix does not depend on the values of the forcing term, it can be multiplied by any arbitrary forcing term to reproduce exactly the discretization error for that forcing term.
The drawback of this independence, however, is that by itself, a force-independent posterior covariance matrix cannot be informative of the force-dependent discretization error.
This poses a problem for Gaussian process--based approaches, since their posterior covariance depends on the observation locations, but not on the observed values at those locations.
In order to obtain a meaningful error estimate, it is necessary to either have a force-dependent prior distribution, or to perform some postprocessing on the posterior covariance matrix to break this independence.
In this work, we have opted for the latter and have shown how by rescaling the eigenvalues of the posterior covariance matrix based on the fine-scale force vector, a distribution can be obtained whose standard deviation corresponds to the discretization error.

One major drawback of the proposed method, as is the case for many probabilistic numerical methods, is its computational cost.
In order to draw samples from the prior and posterior distribution, a fine-scale solve is needed for each sample.
This issue could be solved by approximating or circumventing these fine-scale solves, for example by using iterative solvers with a limited number of iterations, Markov chain Monte Carlo sampling methods or model order reduction.
Furthermore, the formulation in this work has assumed linearity on the partial differential equations and Gaussianity on the prior distribution.
Extensions of the method beyond these assumptions are not trivial.

Finally, the underlying reason for the development of a Bayesian model for finite element discretization error is to allow for the consistent treatment of discretization error through computational pipelines.
In this work, the focus has been on the forward problem, and the fundamentals of our Bayesian formulation of the finite element method.
The demonstration of the method in an inverse modeling or data assimilation context has been left for future work.

\backmatter

\bmhead{Acknowledgements}
This work is supported by the TU Delft AI Labs programme through the SLIMM AI lab.
We are thankful to Uri Peker for the fruitful discussions.

\begin{appendices}
\crefalias{section}{appendix}
\crefalias{subsection}{appendix}

\section{Inhomogeneous boundary conditions}\label{app:inhomogeneous-bcs}
In \cref{subsec:dirichlet-bcs}, only homogeneous Dirichlet boundary conditions were considered.
Here, we expand on this and demonstrate how both homogeneous and inhomogeneous Dirichlet and Neumann boundary conditions can be included in the model.
This is accomplished by assigning a statistically independent normal distribution to the displacement along the Dirichlet boundary $\subd{\vecu}$, as well as the force along the Neumann boundary $\subn{\vecf}$:
\begin{equation}
    \begin{aligned}
        \subd{\vecu} &\sim \normal{\subd{\vecm}}{\beta^2 \subd{\matSigma}} \\
        \subn{\vecf} &\sim \normal{\subn{\vecm}}{\gamma^2 \subn{\matSigma}}
    \end{aligned}
\end{equation}
The assignment of these prior distributions to the Dirichlet and Neumann boundary conditions produces the following prior mean and covariance of the forcing term:
\begin{equation}\label{eq:boundary-prior-moments}
\begin{aligned}
    \boldtilde{\vecm}_\mathbf{f} &= \subid{\matK} \subd{\vecm} + \subn{\vecm} \\
    \boldtilde{\matSigma}_\mathbf{f} &= \matSigmaf + \beta^2 \subid{\matK} \subd{\matSigma} \subid{\matK}^T + \gamma^2 \subn{\matSigma}
\end{aligned}
\end{equation}
By adjusting $\beta$ and $\gamma$, the effects of particular loads can be emphasized or de-emphasized.

This generalization allows for a modeling choice when enforcing inhomogeneous Dirichlet boundary conditions.
These can be strongly enforced in the prior, by setting $\subd{\matSigma} = \vecnull$ and making $\subd{\vecm}$ equal to the true displacement value at the boundary.
Alternatively, they can be weakly enforced by setting $\subd{\vecm} = \vecnull$, and instead assigning a non-zero covariance $\subd{\matSigma}$.
In this case, the Dirichlet boundaries are enforced in a weak sense, because their enforcement is only due to the right-hand side modifications being included in the observations.
Naturally, a combination of these two approaches, where both $\subd{\vecm}$ and $\subd{\matSigma}$ are non-zero is also valid.
For homogeneous Dirichlet boundary conditions, setting $\subd{\vecm} = \vecnull$ and $\subd{\matSigma} = \vecnull$ already strongly enforces the boundary conditions, but this strong enforcement can be weakened by introducing a non-zero $\subd{\matSigma}$.
For Neumann boundary conditions, the same modeling choice between strong and weak enforcement of the boundary conditions applies.

A final point to address is which covariance structure should be applied to the Dirichlet and Neumann covariances $\subd{\matSigma}$ and $\subn{\matSigma}$.
For single point loads and single point constraints, the answer to this question is trivial, namely a null matrix, except for a unit diagonal entry associated with the point load or constraint degree of freedom.
If the problem contains multiple independent point loads or constraints, the covariance structure of $\subd{\matSigma}$ and $\subn{\matSigma}$ is still relatively straightforward:
in this case, the off-diagonal terms of $\subd{\matSigma}$ and $\subn{\matSigma}$ can simply be set to 0.
However, if for example an inhomogeneous Dirichlet or Neumann boundary condition is applied along an edge, this assumption of independence does not hold, and a full covariance structure needs to be obtained for $\subd{\matSigma}$ and $\subn{\matSigma}$.

\section{Sampling the prior and posterior}\label{app:drawing-samples}
The main computational bottleneck of the method lies in the handling of large covariance matrices.
Since both the prior and posterior covariance matrix are full $n \times n$ matrices, their explicit computation, storage and handling quickly becomes infeasible as $n$ increases.
Additionally, $\matK^{-1}$ appears in the expressions for both of these matrices, which suffers from similar problems when computed explicitly.
In this section, methods of sampling exactly from the prior and posterior distributions are discussed.
From these samples, the mean vectors and covariance matrices can be approximated.

\subsection{Ensemble approximation}\label{subapp:ensemble-kalman}
Instead of a mean vector and covariance matrix, an ensemble $\matX$ is used to represent the prior distribution.
$\matX$ is an $n \times N$ matrix containing $N$ samples from the prior distribution.
The prior mean and covariance can be approximated by computing the sample mean $\boldhat{\vecm}$ and sample covariance $\boldhat{\matSigma}$ of the ensemble.
The accuracy of this approximation is depends on the size of the ensemble:
as $N$ increases, the sample mean vector and covariance matrix converge to their exact counterparts, but the computational cost increases accordingly.
By the central limit theorem, the sample mean and covariance will converge to the true sample mean and covariance at a rate of $\frac{1}{\sqrt{N}}$~\citep{berry_accuracy_1941}.

The entries of the sample mean vector and covariance matrix are given by:
\begin{equation}\label{eq:sample-mean}
\begin{aligned}
    \boldhat{\vecm}_i &= \frac{1}{N} \sum_{j=1}^N \matX_{ij} \\
    \boldhat{\matSigma}_{ij} &= \frac{1}{N-1} \sum_{k=1}^N \left(\matX_{ik} - \boldhat{\vecm}_i \right) \left(\matX_{jk} - \boldhat{\vecm}_j \right)
\end{aligned}
\end{equation}
The sample covariance matrix can also be expressed as the following decomposition:
\begin{equation}\label{eq:sample-covariance-decomposition}
\begin{aligned}
    \boldhat{\matSigma} = \frac{1}{N-1} \boldhat{\matF} \boldhat{\matF}^T
\end{aligned}
\end{equation}
Here, $\boldhat{\matF}$ is the sample residual matrix, given by:
\begin{equation}\label{eq:sample-residual-matrix}
\begin{aligned}
    \boldhat{\matF} = \matX - \boldhat{\vecm} \, \vecone_N^T
\end{aligned}
\end{equation}
In the expression above, $\vecone_N$ is a vector of ones of size $N$.

\subsection{Prior sampling}\label{subapp:prior-samples}
One key observation to make about the prior covariance in \cref{eq:noise-prior-covariance} is that the $\matK^{-1}$ term only appears as a premultiplier at the front and as a postmultiplier at the end of the expression.
As a result, instead of obtaining a sample $\boldtilde{\vecu}$ from $\vecu \sim \normal{\vecnull}{\matK^{-1} \matSigmaf \matK^{-1}}$, a sample $\boldtilde{\vecf}$ can be obtained from $\vecf \sim \normal{\vecnull}{\matSigmaf}$.
Since $\matSigmaf$ is a sparse matrix, the Cholesky decomposition of $\matSigmaf$ that is needed to sample $\boldtilde{\vecf}$ is relatively cheap compared to that of a full matrix.
If $\matSigmaf$ is approximated by diagonalizing it, its Cholesky decomposition is trivial.
After obtaining $\boldtilde{\vecf}$, we can compute $\boldtilde{\vecu}$ by solving $\matK \boldtilde{\vecu} = \boldtilde{\vecf}$.

Since multiple solves are needed of the same system, but with a changing right-hand side vector, a direct solver approach is a natural choice.
Here, we use the CHOLMOD library~\citep{davis_user_2013}, which solves sparse linear systems by finding a sparse Cholesky factorization of $\matK$, and then solving for $\boldtilde{\vecu}$ through backward substitution.
Notice that only the cheap backward substitutions need to be repeated for each sample, and the main computational bottleneck of obtaining the factorization needs to be performed only once.
Additionally, for both the white-noise prior from \cref{subsec:white-noise-prior} as well as the Green's function prior from \cref{subsec:greens-function-prior}, $\matSigmaf$ and will have the same sparsity structure as $\matK$.
This means that the permutation that minimizes the fill-ins in the Cholesky decomposition of $\matK$ can be reused for that of $\matSigmaf$.
One drawback of using a direct solver, however, is that it tends to scale poorly as the number of degrees of freedom in the system increases.
For systems with a large number of degrees of freedom, it might therefore be necessary to find alternative sampling methods, though this falls beyond the scope of this paper.

\subsection{\texorpdfstring{Posterior sampling from $\vecu$ via $\vecf$}{Posterior sampling from u via f}}\label{subapp:posterior-samples-from-f}
To obtain posterior samples with the mean and covariance from \cref{eq:noise-posterior-moments}, a similar approach can be applied.
Instead of computing a posterior sample $\boldtilde{\vecu}^*$ directly from $\vecu | \vecg \sim \normal{\boldstar{\vecm}}{\boldstar{\matSigma}}$, a sample of the force vector posterior $\boldtilde{\vecf}^*$ is obtained, from the following distribution
\begin{equation}\label{eq:force-posterior-distribution}
\begin{aligned}
    \vecf | \vecg \sim \normal{\boldstar{\vecmf}}{\boldstar{\matSigmaf}}
\end{aligned}
\end{equation}
where the force posterior mean vector $\boldstar{\vecmf}$ and force posterior covariance matrix $\boldstar{\matSigmaf}$ are given by:
\begin{equation}\label{eq:force-posterior-mean}
\begin{aligned}
    \boldstar{\vecmf} &= \matSigmaf \matPhi \left(\matPhi^T \matSigmaf \matPhi + \sigma_e^2 \matI \right)^{-1} \vecg \\
    &= \matGf \vecg
\end{aligned}
\end{equation}
and
\begin{equation}\label{eq:force-posterior-covariance}
\begin{aligned}
    \boldstar{\matSigmaf} &= \matSigmaf - \matSigmaf \matPhi \left(\matPhi^T \matSigmaf \matPhi + \sigma_e^2 \matI \right)^{-1} \matPhi^T \matSigmaf \\
    &= \left(\matI - \matGf \matPhi^T \right) \matSigmaf
\end{aligned}
\end{equation}
Here, $\matGf$ is the force vector Kalman gain matrix, which is given by:
\begin{equation}\label{eq:force-kalman-gain}
\begin{aligned}
    \matGf = \matSigmaf \matPhi \left(\matPhi^T \matSigmaf \matPhi + \sigma_e^2 \matI \right)^{-1}
\end{aligned}
\end{equation}
If the force posterior covariance matrix is written in its so-called Joseph form, it becomes clear how posterior force vector samples can be computed from prior force vector samples:
\begin{equation}\label{eq:force-posterior-covariance-joseph-form}
\begin{aligned}
    \boldstar{\matSigmaf} &= \left(\matI - \matGf \matPhi^T \right) \matSigmaf \left(\matI - \matGf \matPhi^T \right)^T + \sigma_e^2 \matGf \matGf^T
\end{aligned}
\end{equation}
A prior force vector sample $\boldtilde{\vecf}$, as well as a sample $\boldtilde{\vece}$ from the observation noise $\vece \sim \normal{\vecnull}{\sigma_e^2 \matI}$ are then drawn independently.
The sample from the force vector posterior is then given by:
\begin{equation}\label{eq:force-posterior-sample}
\begin{aligned}
    \boldtilde{\vecf}^* &= \boldstar{\vecmf} + \left(\matI - \matGf \matPhi^T \right) \boldtilde{\vecf} + \matGf \boldtilde{\vece} \\
    &= \boldtilde{\vecf} + \matGf \left(\vecg - \matPhi^T \boldtilde{\vecf} + \boldtilde{\vece} \right)
\end{aligned}
\end{equation}
At this point the sample from the prior distribution of the displacement field can be obtained by solving $\matK \boldtilde{\vecu}^* = \boldtilde{\vecf}^*$.

\subsection{\texorpdfstring{Posterior sampling from $\vecu$ directly}{Posterior sampling from u directly}}\label{subapp:posterior-samples-from-u}
This approach of computing posterior samples by updating prior samples based on perturbed observations of the data $\vecg - \boldtilde{\vece}$ can also be applied to $\vecu$ directly.
Starting from \cref{eq:force-posterior-sample}, $\boldtilde{\vecf}$ and $\matPhi^T$ can be replaced by $\boldtilde{\vecu}$ and $\matH$, respectively.
The force vector Kalman gain matrix given in \cref{eq:force-kalman-gain} can be replaced by the Kalman gain matrix $\matG$ associated with the displacement field:
\begin{equation}\label{eq:kalman-gain}
\begin{aligned}
    \matG &= \matK^{-1} \matSigmaf \matPhi \left(\matPhi^T \matSigmaf \matPhi + \sigma_e^2 \matI \right)^{-1} \\
    &= \matSigma \matH^T \left(\matH \matSigma \matH^T + \sigma_e^2 \matI \right)^{-1}
\end{aligned}
\end{equation}
A prior sample $\boldtilde{\vecu}$ is computed as described in \cref{subapp:prior-samples}, and an independent observation noise sample $\boldtilde{\vece}$ is obtained as before from $\vece \sim \normal{\vecnull}{\sigma_e^2 \matI}$.
The posterior sample $\boldtilde{\vecu}^*$ is then given by:
\begin{equation}\label{eq:posterior-sample}
\begin{aligned}
    \boldtilde{\vecu}^* &= \boldstar{\vecm} + \left(\matI - \matG \matH \right) \boldtilde{\vecu} + \matG \boldtilde{\vece} \\
    &= \boldtilde{\vecu} + \matG \left(\vecg - \matH \boldtilde{\vecu} + \boldtilde{\vece} \right)
\end{aligned}
\end{equation}
Note that this approach requires two $\matK \boldtilde{\vecu} = \boldtilde{\vecf}$ solves to obtain a single sample from the posterior, as opposed to the approach presented in \cref{subapp:posterior-samples-from-f}, where only a single solve is needed.
One reason to still prefer the method of sampling directly from $\vecu$ is that the final solve of $\matK \boldtilde{\vecu}^* = \boldtilde{\vecf}^*$ does not extend well to non-linear problems.
Although non-linear partial differential equations fall beyond the scope of this paper, the option to sample directly from $\vecu$ keeps the door open to this class of problems.

\end{appendices}


\bibliography{citations}

\end{document}